\documentclass[11pt]{article}
\usepackage[utf8]{inputenc}
\usepackage[english]{babel}
\usepackage{graphicx}
\usepackage{amsmath}
\usepackage{amssymb}
\usepackage{amsfonts}
\usepackage{amsthm}
\usepackage[left=2.1cm,top=1.5cm,right=2.1cm, bottom=2.2cm,letterpaper]{geometry}
\usepackage{mathrsfs}
\usepackage{caption}
\usepackage{subcaption}
\usepackage{latexsym}
\usepackage{xfrac}
\usepackage{tcolorbox}
\usepackage{enumitem}
\usepackage{float}
\usepackage{hyperref}
 \usepackage{url}
\usepackage[toc]{appendix}
\newtheorem{theorem}{Theorem}[section]
\newtheorem{corollary}{Corollary}[section]

\newtheorem{definition}{Definition}[section]
\newtheorem{proposition}{Proposition}[section]

\newtheorem{remark}{Remark}[section]
\numberwithin{equation}{section}
\numberwithin{figure}{section}

\newcommand{\R}{\mathbb{R}}

\makeatletter
\renewcommand*\env@matrix[1][*\c@MaxMatrixCols c]{%
	\hskip -\arraycolsep
	\let\@ifnextchar\new@ifnextchar
	\array{#1}}
\makeatother

\def\eq#1{(\ref{#1})}
\def\neweq#1{\begin{equation}\label{#1}}
\def\endeq{\end{equation}}

\begin{document}
	
\title{On the steady motion of Navier-Stokes flows past a fixed obstacle\\
in a three-dimensional channel under mixed boundary conditions}

\author{Gianmarco SPERONE\\
{\small Department of Mathematical Analysis, Charles University in Prague, Czech Republic}}
\date{}
\maketitle
\vspace*{-6mm}
\begin{abstract}
\noindent
We analyze the steady motion of a viscous incompressible fluid in a three-dimensional channel containing an obstacle through the Navier-Stokes equations with mixed boundary conditions: the inflow is given by a fairly general datum and the flow is assumed to satisfy a \textit{constant traction} boundary condition on the outlet, together with the standard no-slip assumption on the obstacle and on the remaining walls of the domain. Explicit bounds on the inflow velocity guaranteeing existence and uniqueness of such steady motion are provided after estimating some Sobolev embedding constants and constructing a suitable solenoidal extension of the inlet velocity through the Bogovskii formula. A quantitative analysis of the forces exterted by the fluid over the obstacle constitutes the main application of our results: by deriving a volume integral formula for the drag and lift, explicit upper bounds on these forces are given in terms of the geometrical constraints of the domain. \par\noindent
\end{abstract}
	
\section{Introduction} \label{motivo}

The steady motion of a viscous incompressible fluid around a three-dimensional fixed obstacle represents an essential and ubiquitous subject-matter in fluid mechanics; in fact, the development of modern aerodynamics was only possible after a proper experimental and theoretical set-up of this problem \cite{ackroyd2001early,bongazspe}. While the standard practice is to mathematically model such motion through the Navier-Stokes equations \cite{ladyzhenskaya1969mathematical}, whether the fluid is allowed to flow in an unbounded region outside the obstacle (the so-called \textit{exterior problem}, see \cite[Chapter XI]{galdi2011introduction} and references therein), or confined to a \textit{bounded} domain, depends on the particular situation considered.
\begin{figure}[b]
	\begin{center}
		\includegraphics[height=55mm,width=80mm]{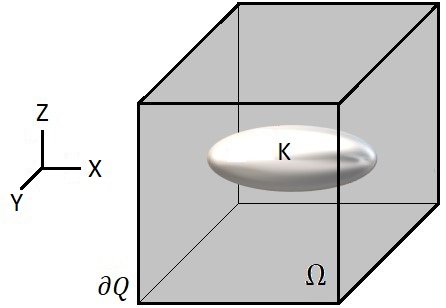}	\qquad \includegraphics[height=55mm,width=80mm]{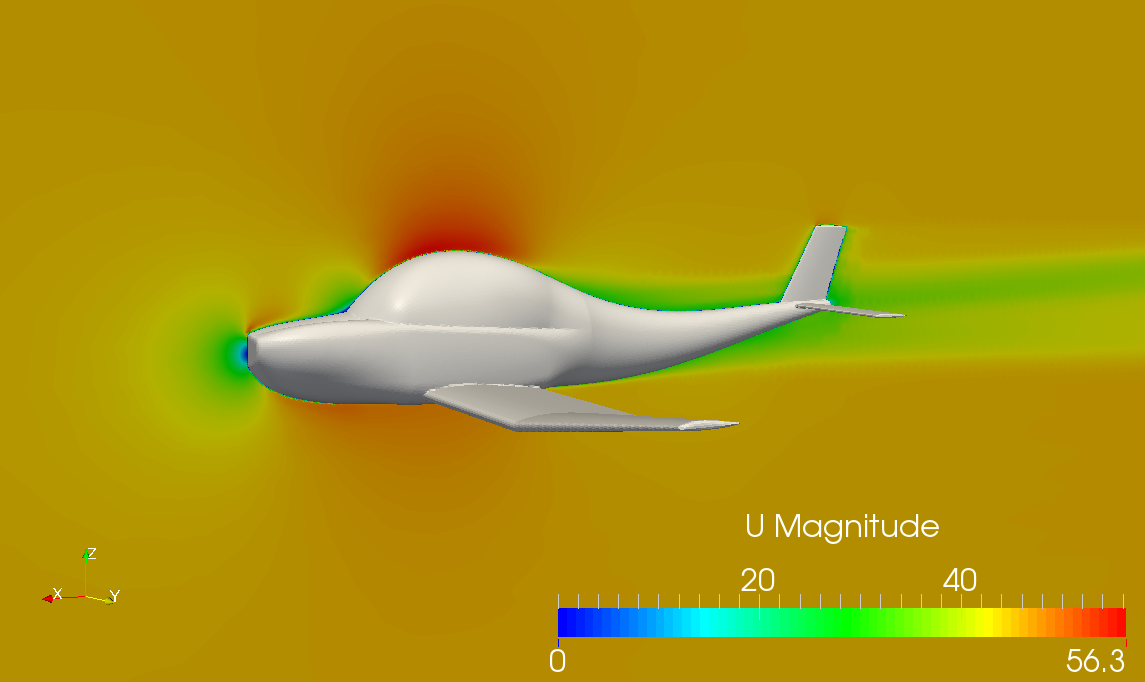}
	\end{center}
	\vspace*{-5mm}
	\caption{Left: the domain $\Omega$ with a smooth obstacle $K$. Right: aerodynamics of a very light aircraft, reproduced by courtesy of Jordi Casacuberta Puig \url{http://www.the-foam-house.es}.}\label{sheet}
\end{figure}
Nevertheless, even if the physical problem naturaly occurs in an unbounded region, the corresponding numerical implementation must be set in a \textit{truncated} bounded domain, therefore creating \textit{artificial boundaries}, see for example \cite{uruba2019reynolds, uruba3d2}. Then, the choice of boundary conditions to be imposed becomes a delicate question, specially on the outlet, because the inlet velocity is usually prescribed and the effects of viscosity dictate that the fluid velocity must equal zero on the surface of the obstacle and on the walls of the domain (no-slip assumption). Since the formulation by Gresho \cite{gresho1991some} in 1991, the \textit{do-nothing} or \textit{constant traction} boundary conditions have been widely employed in Computational Fluid Dynamics (CFD) \cite{braack2014directional, heywood1996artificial, john2002higher, rannacher2012short}, but have also been treated from a strictly mathematical point of view \cite{fursikov2009optimal, kravcmar2001weak, kravcmar2018modeling, kucera1998solutions}, see also the chapter by Galdi in \cite{galdi2008hemodynamical}. Such boundary conditions will be described in detail, see \eqref{nsstokes0}$_4$, after introducing the domain that will be employed hereafter.
\par
In the space $\R^3$ we consider an open, bounded, connected, and simply connected domain $K$, with Lipschitz boundary $\partial K$. Then we
remove $K$, seen as an obstacle, from a larger cube $Q$ such that $\partial K\cap\partial Q=\emptyset$, and we define the domain
\begin{equation} \label{Omega}
Q=(-L,L)^3\, ,\qquad\Omega =Q\setminus \overline{K},
\end{equation}
with $L>\text{diam}(K)$, see Figure \ref{sheet}. We decompose the boundary of $\Omega$ as $ \partial \Omega = \Gamma_{I} \cup \Gamma_{W} \cup \Gamma_{O}$, where
\begin{equation}\label{boundaryomega1}
\begin{aligned}
& \Gamma_{I} = \{-L\} \times [-L,L]^2 \, , \qquad \Gamma_{O} = \{L\} \times [-L,L]^2, \\[3pt]
& \Gamma_{W} = ([-L,L]^2 \times \{\pm L\}) \cup ([-L,L] \times \{\pm L\} \times [-L,L]) \cup \partial K \, .
\end{aligned}
\end{equation}
The outward unit normal to $\partial \Omega$ is denoted by $\hat{n}$. Henceforth we will refer to $\Gamma_{I}$ and $\Gamma_{O}$ in \eqref{boundaryomega1} as the \textit{inlet} and \textit{outlet} of $\Omega$, respectively, while $\Gamma_{W}$ includes all the \textit{physical walls} of $\Omega$.\par
As already mentioned, in the present article we study the steady-state Navier-Stokes equations (in the absence of external forcing) with mixed boundary conditions on the different parts of $\partial \Omega$, that is, the following system of partial differential equations:
\begin{equation}\label{nsstokes0}
\left\{
\begin{aligned}
& -\eta\Delta u+(u\cdot\nabla)u+\nabla p=0,\ \quad  \nabla\cdot u=0 \ \ \mbox{ in } \ \ \Omega, \\[3pt]
& u=h \ \ \mbox{ on } \ \ \Gamma_{I}, \\[3pt]
& u=0 \ \ \mbox{ on } \ \ \Gamma_{W}, \\[3pt]
& -p \hat{n} + \eta \dfrac{\partial u}{\partial n} = 0 \ \ \mbox{ on } \ \ \Gamma_{O}.
\end{aligned}
\right.
\end{equation}
In \eqref{nsstokes0}, $u : \Omega \longrightarrow \mathbb{R}^3$ is the velocity vector field, $p : \Omega \longrightarrow \mathbb{R}$ is the scalar pressure, $\eta>0$ is the (constant) kinematic viscosity and the datum $h : \Gamma_{I} \longrightarrow \mathbb{R}^3$ describes the inlet velocity on $\Gamma_{I}$. For compatibility, $h$ must also vanish on $\partial \Gamma_{I}$. While \eqref{nsstokes0}$_{3}$ describes the usual no-slip boundary condition on the physical walls $\Gamma_{W}$ of $\Omega$, equality \eqref{nsstokes0}$_{4}$ imposes the previously introduced \textit{do-nothing} boundary condition on the outlet. In turns out that such boundary condition arises naturally while studying the weak formulation associated to problem \eqref{nsstokes0}, see Section \ref{boundaryvalueproblems} for more details. Furthermore, from a mechanical point of view, the use of \eqref{nsstokes0}$_{4}$ can be justified because the obstacle perturbs the flow inside the channel, and therefore, it is not possible to know in advance its behavior on the outlet. Nevertheless, as pointed out in \cite{lanzendorfer2020multiple}, boundary conditions such as \eqref{nsstokes0}$_{4}$ create additional difficulties in terms of the well-posedness of \eqref{nsstokes0}; indeed, if $u \in H^{1}(\Omega)$ is a weak solution of \eqref{nsstokes0} (see Definition \eqref{weaksolution}) with $h=0$, integration by parts shows that
\begin{equation} \label{bigproblem}
\eta \, \| \nabla u \|^{2}_{L^{2}(\Omega)} = -\dfrac{1}{2} \int_{\Gamma_{O}} |u|^2 \, (u \cdot e_{1}) \, ,
\end{equation}
where $\{e_{1}, e_{2}, e_{3} \} \subset \R^3$ denotes the standard basis of $\R^3$. Identity \ref{bigproblem} illustrates the fact that the boundary condition \eqref{nsstokes0}$_{4}$ does not allow for a control of the kinetic energy of the fluid flow, as the possibility of a \textit{backwards} flow coming into $\Omega$ from $\Gamma_{O}$ is not excluded. To circumvent this difficulty, problem \eqref{nsstokes0} with additional assumptions on the normal component $u \cdot \hat{n}$ on $\Gamma_{O}$ has been studied in \cite{braack2014directional, bruneau1996new, feistauer2013existence, kravcmar2001weak}. Otherwise, existence and uniqueness of solutions for \eqref{nsstokes0} can be proved only for small data; examples of non-existence or multiplicity of solutions were provided in \cite{lanzendorfer2020multiple}. More precisely, as explained in \cite[Chapter III]{galdi2008hemodynamical}, smallness of the inlet velocity guarantees existence of a unique solution \textit{within some ball}, see also Theorem \ref{mainresult}. It remains an open problem to determine whether a given solution of \eqref{nsstokes0}, inside a given ball, is unique in the class of \textit{all} possible weak solutions corresponding to the same data.
\par
The present articles serves as a contribution in the theoretical analysis of \eqref{nsstokes0}, with the purpose of conferring a deeper understanding of some of the questions described before and motivating further studies on this matter. The main goal of this work regards the \textit{quantification} of the smallness assumption on the inlet velocity that ensures the well-posedness of problem \eqref{nsstokes0}: an explicit upper bound on the \textit{size} of the boundary datum $h$, in terms of the geometric constraints of the obstacle $K$ and of the domain $\Omega$, is given in Section \ref{boundaryvalueproblems}, see Theorem \ref{mainresult} and Corollary \ref{maincor}, which then also yield an explicit upper bound for the solution of \eqref{nsstokes0}. Similar results were obtained in \cite{filippoclara,gazspe} regarding the unique solvability of \eqref{nsstokes0}$_1$ under non-homogeneous Dirichlet boundary conditions on $\partial Q$. By following closely the proofs given by Fursikov \& Rannacher in \cite{fursikov2009optimal}, one discovers that, in order to yield such explicit threshold, several tools of functional analysis must be carefully studied. Therefore, in Section \eqref{secsobolev} we firstly provide lower bounds for the Sobolev constant of the embedding $H^{1}(\Omega) \subset L^{p}(\Omega)$, $p \in [2,6]$, involving functions that vanish only on $\Gamma_{W}$ or on $\Gamma_{I} \cup \Gamma_{W}$. Thus, symmetrization techniques as in \cite{filippoclara, cortona, gazspe} can only be applied after properly reflecting $\Omega$ with respect to the planes $x= \pm L$, and performing a suitable even extension of the functions considered, as in \cite[Chapter 2]{uraltseva2014linear}. Secondly, in Section \ref{bogovskii} we build an \textit{estimable} solenoidal extension of the boundary datum $h$, namely, we look for a vector field $v_{0} \in H^{1}(\Omega)$ such that 
\begin{equation}\label{div}
\nabla\cdot v_0=0 \ \text{ in } \ \Omega, \qquad v_0=h \ \text{ on } \ \Gamma_{I}, \qquad v_0=0 \ \text{ on } \ \Gamma_{W}.
\end{equation}
Problem \eqref{div} has been extensively studied in the past, not only because of its applicability in fluid mechanics (see the works by Ladyzhenskaya \& Solonnikov \cite{ladyzhenskaya1969mathematical,ladyzhenskaya1978some}, Bogovskii \cite{bogovskii1979solution} and the book by Galdi \cite[Section III.3]{galdi2011introduction}), but also due to its purely mathematical interest and connection with the Calderón-
Zygmund theory of singular integrals, see the book by Acosta \& Durán \cite{acosta2017divergence}.
Our construction invokes the method described in \cite{gazspefra}: by inverting the trace operator, an extension of $h$ (not necessarily solenoidal) is determined; then the Bogovskii problem with the resulting divergence is studied to obtain a solenoidal extension; finally, by solving a variational problem involving the infinity-Laplacian and using ad hoc cut-off functions, an explicit estimate of $v_{0}$ is given.\par 
As an application of our results, in Section \ref{forces} we study the fluid forces exerted over $K$. The \textit{lift force}, understood as the force component that is perpendicular to the oncoming stream, plays a fundamental role in aerodynamics (where it must be maximized in order counter the force of gravity acting over the aircraft \cite[Chapter 3]{ackroyd2001early}) and in civil engineering (where it must be minimized in order to avoid instabilities of structures such as suspension bridges or skyscrapers \cite{gazzola2015mathematical}). Therefore, it would be desirable to have an explicit upper bound on the lift, in terms of the known data; from a mathematical point of view, this represents a difficult task, due to the lack of estimates regarding the continuity constant of the trace operator. Nevertheless, by building an auxiliary solenoidal extension of constant boundary data on $\partial K$ (namely, of the unit vectors $e_1$ and $e_3$, see Proposition \ref{solext}), in Proposition \ref{newformula} we derive a volume integral formula for both the drag and lift; after combining such formula with the results of Section \ref{boundaryvalueproblems}, explicit upper bounds for the drag and lift are yielded. In order to highlight the utility of this result, a concrete example is developed in which the inlet velocity is given by an analytic expression, see Corollary \ref{upperdraglift}.\par
This paper is organized as follows. In Section \ref{funcine} we exploit some techniques of functional analysis that are needed for the study of \eqref{nsstokes0}, such as: inequalities associated to the immersion $H^1(\Omega)\subset L^p(\Omega)$, for $p \in [2,6]$ (regarding functions that vanish on a proper part of $\partial \Omega$), and a priori bounds for the solutions of \eq{div}. In Section \ref{boundaryvalueproblems} we begin by providing estimates for the (linear) Stokes operator under mixed boundary conditions in the presence of an external force, which are essential to establish the well-posedness of \eqref{nsstokes0} and to provide the explicit estimates given in the main results (Theorem \ref{mainresult} and Corollary \ref{maincor}). Finally, Section \ref{forces} is devoted to a quantitative study of the drag and lift forces applied by the fluid over $K$, the main goal being the obtainment of explicit upper bounds for such forces in terms of the geometric parameters of $K$.

\section{Tools from functional analysis} \label{funcine}
We emphasize that, for the sake of simplicity, no distinction will be made for the notation of functional spaces of scalars, vectors or matrices.
\subsection{Explicit bounds for some Sobolev embedding constants} \label{secsobolev}
Let $\Omega$ be as in \eqref{Omega}. We consider the following Sobolev spaces of functions vanishing on subsets of $\partial \Omega$ having a positive 2D-measure:
$$
V_{*}(\Omega) = \{ v \in H^{1}(\Omega) \ | \ v=0 \ \ \mbox{on} \ \ \Gamma_{I} \cup \Gamma_{W} \} \, , \qquad V(\Omega) = \{ v \in H^{1}(\Omega) \ | \ v=0 \ \ \mbox{on} \ \ \Gamma_{W} \};
$$
notice that $V_{*}(\Omega) \subset V(\Omega)$. Since $|\Gamma_{W}|>0$, the Poincaré inequality holds in $V(\Omega)$, which means that $v\mapsto\|\nabla v\|_{L^2(\Omega)}$ is indeed a norm on $V(\Omega)$ (and also on $V_{*}(\Omega)$). Moreover, in view of the embedding $H^{1}(\Omega) \subset L^{p}(\Omega)$ for every $p \in [2,6]$, the Sobolev constants for the embeddings $V_{*}(\Omega) \subset L^{p}(\Omega)$ and $V(\Omega) \subset L^{p}(\Omega)$ admit a variational definition:
\begin{equation} \label{sobolevconstants1}
\mathcal{S}_{p} \doteq \min_{w\in V_{*}(\Omega) \setminus \{0\}} \ \dfrac{\|\nabla w\|^2_{L^2(\Omega)}}{\|w\|^{2}_{L^{p}(\Omega)}} \qquad \text{and} \qquad \mathcal{J}_{p} \doteq \min_{w\in V(\Omega) \setminus \{0\}} \ \dfrac{\|\nabla w\|^2_{L^2(\Omega)}}{\|w\|^{2}_{L^{p}(\Omega)}} \qquad \forall p \in [2,6].
\end{equation}
For every $p \in [2,6]$ we have $0 < \mathcal{J}_{p} \leq \mathcal{S}_{p}$ and, by \eqref{sobolevconstants1},
\begin{equation} \label{sobolevconstants11}
\mathcal{S}_{p} \, \|v\|^{2}_{L^{p}(\Omega)} \leq \|\nabla v\|^2_{L^2(\Omega)} \qquad \forall v \in V_{*}(\Omega); \qquad \mathcal{J}_{p} \, \|v\|^{2}_{L^{p}(\Omega)} \leq \|\nabla v\|^2_{L^2(\Omega)} \qquad \forall v \in V(\Omega).
\end{equation}
We emphasize that the inequalities in \eqref{sobolevconstants11} are equally valid for scalar or vector functions (with the same constant). Indeed, if $v = (v_1, v_2, v_3) \in V(\Omega)$ is a vector field, by the Minkowski inequality we get
$$
\begin{aligned}
\|v\|^{p}_{L^p(\Omega)} & = \left\| \, | v_{1} |^{2} + | v_{2} |^{2} + | v_{3} |^{2} \, \right\|^{p/2}_{L^{p/2}(\Omega)} \leq \left( \|v_{1}\|^{2}_{L^p(\Omega)} + \|v_{2}\|^{2}_{L^p(\Omega)} + \|v_{3}\|^{2}_{L^p(\Omega)} \right)^{p/2} \\[3pt]
& \leq \left( \dfrac{1}{\mathcal{J}_{p}} \right)^{p/2} \left( \| \nabla v_{1}\|^{2}_{L^2(\Omega)} + \| \nabla v_{2}\|^{2}_{L^2(\Omega)} + \| \nabla v_{3}\|^{2}_{L^2(\Omega)} \right)^{p/2} = \left( \dfrac{1}{\mathcal{J}_{p}} \right)^{p/2} \| \nabla v \|^{p}_{L^2(\Omega)} \, .
\end{aligned}
$$
The purpose of this section is to provide explicit lower bounds for the constants $\mathcal{S}_{p}$ and $\mathcal{J}_{p}$, for $p \in [2,6]$. We firstly treat the case $p=2$ (Poincaré inequality) and the limiting case $p=6$, where the value for the best Sobolev constant of the embedding $H_{0}^{1}(\Omega) \subset L^{6}(\Omega)$ is well-known and does not depend on the domain considered, see \cite{lions1988best, talenti1976best}.
 
\begin{theorem}\label{boundsobolev}
Let $\Omega$ be as in \eq{Omega}. For any scalar or vector function $u\in V_{*}(\Omega)$ one has
\begin{equation} \label{ineqL4H1u2}
\|u\|_{L^2(\Omega)} \leq \dfrac{1}{\pi} \min \left\{ \sqrt[3]{\frac{3}{2\pi} (|Q|-|K|)}  \, , \,  \dfrac{4L}{3} \right\}  \, \|\nabla u\|_{L^2(\Omega)} \, , \qquad \|u\|_{L^6(\Omega)} \leq \dfrac{2}{\sqrt{3} \, \pi^{2/3}} \, \|\nabla u\|_{L^2(\Omega)} \, .
\end{equation}
For any scalar or vector function $u\in V(\Omega)$ one has
\begin{equation} \label{ineqL4H1u22}
\begin{aligned}
\|u\|_{L^2(\Omega)} & \leq \dfrac{\sqrt{3}}{\left( \pi^2 \max \left\{ \sqrt[3]{\dfrac{2 \pi}{3 (|Q|-|K|)} }  \, , \,  \dfrac{3}{4L} \right\}^{2} - \dfrac{2}{L^2} \right)^{1/2}} \, \|\nabla u\|_{L^2(\Omega)} \, , \\[8pt] 
\|u\|_{L^6(\Omega)} & \leq  \dfrac{ 2 \pi^{1/3}  \max \left\{ \sqrt[3]{\dfrac{2 \pi}{3 (|Q|-|K|)} }  \, , \,  \dfrac{3}{4L} \right\}}{\left( \pi^2 \max \left\{ \sqrt[3]{\dfrac{2 \pi}{3 (|Q|-|K|)} }  \, , \,  \dfrac{3}{4L} \right\}^{2} - \dfrac{2}{L^2} \right)^{1/2}} \, \|\nabla u\|_{L^2(\Omega)} \, .
\end{aligned}
\end{equation}
\end{theorem}
\noindent
\begin{proof}
Let us denote by $\Omega_{\sharp} \subset \mathbb{R}^3$ the reflection of $\Omega$ with respect to the plane $x=L$, that is,
\begin{equation} \label{reflectomega}
\Omega_{\sharp} = \{ (2L-x,y,z) \ | \ (x,y,z) \in \Omega \}.
\end{equation}
In the same way $K_{\sharp} \subset  Q_{\sharp} \subset \mathbb{R}^3$ are defined as the reflections of $K$ and $Q$, correspondingly, with respect to the plane $x=L$:
\begin{equation} \label{reflectomega2}
K_{\sharp} = \{ (2L-x,y,z) \ | \ (x,y,z) \in K \} \, , \qquad Q_{\sharp} = (L,3L) \times (-L,L)^2 \, .
\end{equation}
We then write $\Omega_{+} = \Omega \cup \Omega_{\sharp}$ and, given a function $u \in V_{*}(\Omega)$ (scalar or vector), we define its even extension to $\Omega_{+}$ according to the formula
\begin{equation} \label{evenext1}
u_{+}(x,y,z) =
\left\{
\begin{array}{ll}
u(x,y,z) \ \ & \text{if} \ (x,y,z) \in \Omega \\[3pt]
u(2L-x,y,z) \ \ & \text{if} \ (x,y,z) \in \Omega_{\sharp} \, ,
\end{array}
\right.
\end{equation}
so that $u_{+} \in H^{1}_{0}(\Omega_{+})$. Through the Faber-Krahn inequality \cite{faber1923beweis, krahn1925rayleigh} we first bound the $L^2(\Omega_{+})$-norm of $u_{+}$ in terms of its Dirichlet norm by using the Poincaré inequality in $\Omega^*$, namely a ball having the same measure as $\Omega_{+}$. Since $|\Omega_{+}|=2(|Q|-|K|)$,
the radius of $\Omega^*$ is given by
$$
R=\sqrt[3]{\dfrac{3 |\Omega_{+}|}{4 \pi}}=\sqrt[3]{\frac{3}{2\pi} (|Q|-|K|)} \, .
$$
Since the Poincaré constant (least eigenvalue of $-\Delta$) in the unit ball is given by $\pi$ (corresponding to the first zero of the spherical Bessel
function of order 0), the Poincaré constant of $\Omega^*$ is given by $\pi^2/R^2$, which means that
$$
\min_{w\in H^1_0(\Omega_{+})}\ \frac{\|\nabla w\|_{L^2(\Omega_{+})}}{\|w\|_{L^2(\Omega_{+})}}\, \ge\, \min_{w\in H^1_0(\Omega^*)}\ \frac{\|\nabla w\|_{L^{2}(\Omega^{*})}}{\|w\|_{L^{2}(\Omega^{*})}}
=\frac{\pi}{R}\, .
$$
Therefore,
\begin{equation} \label{poinc1}
\|u_{+}\|_{L^2(\Omega_{+})} \leq \frac{R}{\pi}\|\nabla u_{+}\|_{L^2(\Omega_{+})}=\frac{1}{\pi}\, \sqrt[3]{\frac{3}{2\pi} (|Q|-|K|)}  \, \|\nabla u_{+}\|_{L^2(\Omega_{+})} \, .
\end{equation}
On the other hand, the function
$$
\cos \left( \dfrac{\pi}{4L}(x-L) \right) \cos \left( \dfrac{\pi y}{2L} \right) \cos \left( \dfrac{\pi z}{2L} \right) \qquad \forall (x,y,z) \in Q_{+} \doteq (-L,3L) \times (-L,L)^2
$$
is a non-negative eigenfunction of the operator $-\Delta$ in $Q_{+}$ under homogeneous Dirichlet boundary conditions. As a consequence, it is associated to the least eigenvalue which is then given by $\lambda = \frac{9 \pi^{2}}{16 L^{2}}$. Therefore, the Poincar\'e inequality in $Q_{+}$ reads
$$
\|w\|_{L^2(Q_{+})} \leq \frac{4 L}{3 \pi} \, \|\nabla w\|_{L^2(Q_{+})} \qquad \forall w\in H^1_0(Q_{+}) \, .
$$
Since $u_{+}$ can be extended by 0 in $K$ and in $K_{\sharp}$, it becomes an element of $H^1_0(Q_{+})$ that satisfies
\begin{equation} \label{poinc2}
\|u_{+} \|_{L^2(\Omega_{+})} \leq \frac{4 L}{3 \pi} \, \|\nabla u_{+} \|_{L^2(\Omega_{+})} \, .
\end{equation}
Now, it can be easily seen that
\begin{equation} \label{poinc3}
\|u_{+} \|^{2}_{L^2(\Omega_{+})} = 2 \, \|u \|^{2}_{L^2(\Omega)} \, , \qquad  \|u_{+} \|^{6}_{L^6(\Omega_{+})} = 2 \, \|u \|^{6}_{L^6(\Omega)} \, , \qquad \|\nabla u_{+} \|^{2}_{L^2(\Omega_{+})} = 2 \, \|\nabla u \|^{2}_{L^2(\Omega)}  \, ,
\end{equation}
which, inserted respectively into \eqref{poinc1} and \eqref{poinc2}, yields the two bounds in \eqref{ineqL4H1u2}$_1$.\par
Regarding the critical exponent $p=6$, from the work of Talenti \cite{talenti1976best} we obtain the optimal Sobolev constant for the embedding $H_{0}^{1}(\Omega_{+}) \subset L^{6}(\Omega_{+})$, so that
\begin{equation} \label{talenti1}
\|u_{+}\|_{L^6(\Omega_{+})} \leq \dfrac{4^{1/3}}{\sqrt{3} \, \pi^{2/3}} \, \|\nabla u_{+}\|_{L^2(\Omega_{+})} \, .
\end{equation}
Inequality \eqref{ineqL4H1u2}$_2$ is then obtained by inserting \eqref{poinc3} into \eqref{talenti1}.\par 
Consider now a scalar function $u \in V(\Omega)$, and define $U_{1}, U_{2} \in H^{1}(\Omega)$ according to
$$
U_{1}(x,y,z) = \dfrac{L+x}{2L} u(x,y,z) \, , \qquad U_{2}(x,y,z) = \dfrac{L-x}{2L} u(x,y,z) \qquad \text{for a.e.} \ \ (x,y,z) \in \Omega \, .
$$
Since $U_{1}$ vanishes on $\Gamma_{I} \cup \Gamma_{W}$ and $U_{2}$ vanishes on $\Gamma_{O} \cup \Gamma_{W}$, the previous argument (of reflecting $\Omega$ with respect to the planes $x=\pm L$) can be applied to deduce that both $U_{1}$ and $U_{2}$ satisfy the inequalities in \eqref{ineqL4H1u2}. Moreover, since $u= U_{1} + U_{2}$ in $\Omega$, we then obtain
$$
\|u\|^{2}_{L^2(\Omega)} = \|U_{1}\|^{2}_{L^2(\Omega)} + \|U_{2}\|^{2}_{L^2(\Omega)} + \dfrac{1}{2 L^{2}} \int_{\Omega} (L^{2} - x^{2}) |u|^2 \leq \|U_{1}\|^{2}_{L^2(\Omega)} + \|U_{2}\|^{2}_{L^2(\Omega)} + \dfrac{1}{2} \|u\|^{2}_{L^2(\Omega)} \, ,
$$ 
so that
\begin{equation} \label{sperone1}
\begin{aligned}
\dfrac{1}{2} \|u\|^{2}_{L^2(\Omega)} & \leq \|U_{1}\|^{2}_{L^2(\Omega)} + \|U_{2}\|^{2}_{L^2(\Omega)} \\[6pt]
& \leq \left( \dfrac{1}{\pi} \right)^{2} \min \left\{ \sqrt[3]{\frac{3}{2\pi} (|Q|-|K|)}  \, , \,  \dfrac{4L}{3} \right\}^{2} \left( \| \nabla U_{1}\|^{2}_{L^2(\Omega)} + \| \nabla U_{2}\|^{2}_{L^2(\Omega)} \right).
\end{aligned}
\end{equation}
On the other hand, after applying Young's inequality we get
\begin{equation} \label{sperone2}
\begin{aligned}
\| \nabla U_{1}\|^{2}_{L^2(\Omega)} + \| \nabla U_{2}\|^{2}_{L^2(\Omega)} & = \dfrac{1}{2 L^{2}} \|u\|^{2}_{L^2(\Omega)} + \dfrac{1}{L^2} \int_{\Omega} x  u \dfrac{\partial u}{\partial x} + \dfrac{1}{2 L^{2}} \int_{\Omega} (L^{2} + x^{2}) |\nabla u|^2 \\[6pt]
& \leq \dfrac{1}{2 L^{2}} \|u\|^{2}_{L^2(\Omega)} + \dfrac{1}{L} \left( \dfrac{1}{2 L} \|u\|^{2}_{L^2(\Omega)} + \dfrac{L}{2} \| \nabla u\|^{2}_{L^2(\Omega)} \right) + \| \nabla u \|^{2}_{L^2(\Omega)} \\[6pt]
& =  \dfrac{1}{L^{2}} \|u\|^{2}_{L^2(\Omega)} + \dfrac{3}{2} \| \nabla u\|^{2}_{L^2(\Omega)} \, ,
\end{aligned}
\end{equation}
which, once inserted into \eqref{sperone1}, yields \eqref{ineqL4H1u22}$_1$ after noticing that
$$
\left( \dfrac{1}{\pi} \right)^{2} \min \left\{ \sqrt[3]{\frac{3}{2\pi} (|Q|-|K|)}  \, , \,  \dfrac{4L}{3} \right\}^{2} < \dfrac{L^2}{2} \, .
$$
In order to prove \eqref{ineqL4H1u22}$_2$, we start by applying Minkowski's and H\"older's inequality in the following way:
$$
\begin{aligned}
\|u\|^{6}_{L^{6}(\Omega)} & = \left\| \, | U_{1} |^{2} + 2 U_{1} U_{2} + | U_{2} |^{2} \, \right\|^{3}_{L^{3}(\Omega)} \leq \left( \| U_{1} \|^{2}_{L^{6}(\Omega)} + 2 \| U_{1} U_{2} \|_{L^{3}(\Omega)} + \| U_{2} \|^{2}_{L^{6}(\Omega)} \right)^{3} \\[6pt]
& \leq \left( \| U_{1} \|^{2}_{L^{6}(\Omega)} + 2 \| U_{1} \|_{L^{6}(\Omega)} \| U_{2} \|_{L^{6}(\Omega)} + \| U_{2} \|^{2}_{L^{6}(\Omega)} \right)^{3} \\[6pt]
& \leq 8 \left( \| U_{1} \|^{2}_{L^{6}(\Omega)} + \| U_{2} \|^{2}_{L^{6}(\Omega)} \right)^{3} \, .
\end{aligned}
$$
Therefore, in view of \eqref{ineqL4H1u2}$_2$ and \eqref{sperone2} we have
$$
\|u\|_{L^{6}(\Omega)} \leq \dfrac{2 \sqrt{2}}{\sqrt{3} \, \pi^{2/3}} \left( \| \nabla U_{1} \|^{2}_{L^{2}(\Omega)} + \| \nabla U_{2} \|^{2}_{L^{2}(\Omega)} \right)^{1/2} \leq \dfrac{2 \sqrt{2}}{\sqrt{3} \, \pi^{2/3}} \left( \dfrac{1}{L^{2}} \|u\|^{2}_{L^2(\Omega)} + \dfrac{3}{2} \| \nabla u\|^{2}_{L^2(\Omega)} \right)^{1/2} \, ,
$$
from where we derive inequality \eqref{ineqL4H1u22}$_2$ after applying \eqref{ineqL4H1u22}$_1$ in the right-hand side.
\end{proof}

\begin{remark}
The break-even in the bounds \eqref{ineqL4H1u2}$_1$ and \eqref{ineqL4H1u22} occurs whenever
$$
\dfrac{|Q|}{|Q| - |K|} = \dfrac{81}{16 \pi},
$$
and therefore, it is the relative size of $K$ within $Q$ that determines the optimal bound in \eqref{ineqL4H1u2}$_1$ and \eqref{ineqL4H1u22}.
\end{remark}

\begin{remark}
Since $\Omega \subset \{ (x,y,z) \in \mathbb{R}^{3} \ | -L \leq z \leq L \, \}$, the standard proof of the Poincaré inequality for functions in $H^1_0(\Omega)$ (see, for example, \cite[Section II.5]{galdi2011introduction}) yields 
$$
\|w \|_{L^2(\Omega)} \leq L  \, \|\nabla w \|_{L^2(\Omega)} \qquad \forall w \in H^1_0(\Omega) \, ,
$$
and thus, a larger upper bound than the ones given in \eqref{ineqL4H1u2}$_1$-\eqref{ineqL4H1u22}$_1$. In fact, from \eqref{ineqL4H1u2}$_1$- \eqref{ineqL4H1u22}$_1$ we infer
$$
\|v\|_{L^{2}(\Omega)} \leq 0.43 L \|\nabla v\|_{L^2(\Omega)} \qquad \forall v \in V_{*}(\Omega); \qquad  \|v\|_{L^{2}(\Omega)} \leq  0.46 L \|\nabla v\|_{L^2(\Omega)} \qquad \forall v \in V(\Omega).
$$
\end{remark}

Given any $p \in (2,6]$, for later use we define
\begin{equation} \label{delpinoconstant}
C(p) = \left( \dfrac{p^{2}-4}{24 \pi} \right)^{\frac{3(p-2)}{4p}} \left[ \dfrac{\Gamma\left( \dfrac{p+2}{p-2} \right)}{\Gamma\left( \dfrac{10-p}{2(p-2)} \right)} \right]^{\frac{p-2}{2p}} \left[ \dfrac{10-p}{2(p+2)} \right]^{\frac{10-p}{4p}} \, .
\end{equation}
We are now in position to prove the following:
\begin{theorem}\label{boundsobolev2}
Let $\Omega$ be as in \eq{Omega} and $p \in (2,6]$. For any scalar or vector function $u\in V_{*}(\Omega)$ one has
\begin{equation} \label{sobolevine1}
\|u\|_{L^p(\Omega)} \leq 2^{\frac{p-2}{2p}} C(p) \, \left( \dfrac{1}{\pi}  \, \min \left\{ \sqrt[3]{\frac{3}{2\pi} (|Q|-|K|)}  \, , \,  \dfrac{4L}{3} \right\} \right)^{\frac{4p - p^{2} + 12}{2p(p+2)}} \, \|\nabla u\|_{L^2(\Omega)} \, ,
\end{equation}
where $C(p)>0$ is defined in \eqref{delpinoconstant}.
\end{theorem}
\noindent
\begin{proof}
As in the proof of Theorem \ref{boundsobolev}, we introduce the reflected domains (with respect to the plane $x=L$) $\Omega_{\sharp} \subset \mathbb{R}^3$ and $K_{\sharp} \subset  Q_{\sharp} \subset \mathbb{R}^3$ given by \eqref{reflectomega} and \eqref{reflectomega2}. We then write $\Omega_{+} = \Omega \cup \Omega_{\sharp}$ and, given a scalar function $u \in V_*(\Omega)$, we define its (even) extension $u_{+} \in H^{1}_{0}(\Omega_{+})$ according to \eqref{evenext1}. Next, we recall that del Pino-Dolbeault \cite[Theorem 1]{delpino} obtained the following (optimal) Gagliardo-Nirenberg inequality in $\R^3$:
\begin{equation} \label{gagliardo0}
\|w\|_{L^p(\Omega_{+})} \leq C(p)^{\frac{4}{10-p}} \, \|\nabla w\|_{L^2(\Omega_{+})}^{\frac{6(p-2)}{p(10-p)}} \, \|w\|_{L^{\frac{p}{2}+1}(\Omega_{+})}^{\frac{4p - p^{2} + 12}{p(10-p)}} \qquad \forall w \in H^1_0(\Omega_{+}).
\end{equation}
Since functions in $H^1_0(\Omega_{+})$ may be extended by zero outside $\Omega_{+}$, they can be seen as functions defined over the whole space. Therefore, $u_{+}$ also verifies \eqref{gagliardo0}. An application of H\"older's inequality allows us to get rid of the $L^{\frac{p}{2}+1}(\Omega_{+})$-norm:
$$
\| u_{+} \|_{L^{\frac{p}{2}+1}(\Omega_{+})}^{\frac{p}{2}+1} = \int_{\Omega_{+}} |u_{+}|^{\frac{p}{2}} |u_{+}| \leq \| u_{+} \|_{L^{p}(\Omega_{+})}^{\frac{p}{2}} \, \| u_{+} \|_{L^{2}(\Omega_{+})} \,
$$
which, inserted into \eqref{gagliardo0} (replacing $w$ by $u_{+}$), yields
$$
\| u_{+} \|_{L^p(\Omega_{+})} \leq C(p)^{\frac{4}{10-p}} \, \|\nabla u_{+} \|_{L^2(\Omega_{+})}^{\frac{6(p-2)}{p(10-p)}} \, \| u_{+} \|_{L^{p}(\Omega_{+})}^{\frac{4p - p^{2} + 12}{(p+2)(10-p)}} \, \| u_{+} \|_{L^{2}(\Omega_{+})}^{\frac{2(4p - p^{2} + 12)}{p(p+2)(10-p)}} \, ,
$$
or equivalently
\begin{equation} \label{gagliardo1}
\| u_{+} \|^{\frac{4}{10-p}}_{L^p(\Omega_{+})} \leq C(p)^{\frac{4}{10-p}} \, \|\nabla u_{+} \|_{L^2(\Omega_{+})}^{\frac{6(p-2)}{p(10-p)}} \, \| u_{+} \|_{L^{2}(\Omega_{+})}^{\frac{2(4p - p^{2} + 12)}{p(p+2)(10-p)}} \, .
\end{equation}
As in \eqref{poinc3}, it can be easily seen that $\|u_{+} \|^{p}_{L^p(\Omega_{+})} = 2 \, \|u \|^{p}_{L^p(\Omega)}$, so that \eqref{gagliardo1} becomes
\begin{equation} \label{gagliardo2}
\| u \|^{\frac{4}{10-p}}_{L^p(\Omega)} \leq 4^{\frac{p-2}{p(10-p)}} \, C(p)^{\frac{4}{10-p}} \, \|\nabla u \|_{L^2(\Omega)}^{\frac{6(p-2)}{p(10-p)}} \, \| u \|_{L^{2}(\Omega)}^{\frac{2(4p - p^{2} + 12)}{p(p+2)(10-p)}} \, .
\end{equation}
Inequality \eqref{sobolevine1} follows by inserting \eqref{ineqL4H1u2}$_1$ into \eqref{gagliardo2}, and then taking the $\frac{4}{10-p}$-roots of the resulting inequality.
\end{proof}

As a direct consequence of Theorem \ref{boundsobolev2} we obtain
\begin{corollary}\label{boundsobolev3}
Let $\Omega$ be as in \eq{Omega}. For any scalar or vector function $u \in V_{*}(\Omega)$ one has
\begin{equation} \label{sobolevl3}
\|u\|_{L^3(\Omega)} \leq \dfrac{1}{\sqrt{5}} \left( \dfrac{7^{7/3}}{3} \right)^{1/4} \left( \dfrac{1}{2 \pi^{5}} \right)^{1/6} \, \min \left\{ \sqrt[3]{\frac{3}{2\pi} (|Q|-|K|)}  \, , \,  \dfrac{4L}{3} \right\}^{1/2} \, \|\nabla u\|_{L^2(\Omega)} \, .
\end{equation}
\end{corollary}

\begin{remark}
Notice that, by putting $p=6$ in \eqref{sobolevine1}, we also recover inequality \eqref{ineqL4H1u2}$_2$.
\end{remark}

\subsection{Solenoidal extension of the inlet velocity through Bogovskii's formula} \label{bogovskii}
As we will see in Section \ref{boundaryvalueproblems}, the well-posedness of problem \eqref{nsstokes0} relies on the use of a \textit{solenoidal extension} of the boundary datum $h$, namely, on the existence of a vector field $v_{0} \in H^{1}(\Omega)$ satisfying \eqref{div}. The purpose of this section is to analyze the solvability of \eqref{div} and to obtain {\it explicit bounds} of its solution in terms of the known data (that is, in terms of $h$ and of the geometric parameters of $\Omega$). We follow closely the procedure described in \cite{gazspefra}: we fix real numbers $a, b, c$ such that
\begin{equation}\label{OmegaR}
L >  a \geq b \geq c >0  \quad \text{ and } \quad K \subset  P  \subset Q \, , \quad \text{ with } \ P=(-a,a)\times(-b,b) \times (-c, c) \, ,
\end{equation}so that the obstacle $K$ is enclosed by the parallelepiped $P$. We then set
{\small
\begin{equation}\label{sigmagamma3d}
\begin{aligned}
\sigma & =    7 L^3 -  (a + b + c) L^2 - (a b  +a c + b c) L -a b c , \quad \gamma =   6 L^3- 2 (a + b + c) L^2 - 2 (a b + a c + b c) L+6 a b c , \\[6pt]
 M & = \sqrt{ 12  \left( 1 + \dfrac{16}{\gamma} (L^3 - abc) \right) } \Bigg[ 327.23 + \dfrac{445.17 \sqrt{\sigma}}{(L-a)^{3/2}} + \dfrac{153.85 \sigma}{(L-a)^3}
+ \dfrac{144 L^{2}}{(L-a)^2} \left( 22.4 + \dfrac{15.79 \sqrt{\sigma}}{(L-a)^{3/2}}  \right)^{2} \Bigg]^{1/2},
\end{aligned}
\end{equation}}and may now prove the following:
\begin{theorem} \label{teononconstant}
Assume \eqref{OmegaR} and define $M>0$ as in \eqref{sigmagamma3d}. Then, given $h = (h_1, h_2, h_3) \in H^{1}_{0}(\Gamma_{I})$, there exists a vector field $v_0\in H^{1}(\Omega)$ satisfying \eqref{div} such that
\begin{equation}\label{f:gb}
\| \nabla v_{0} \|_{L^{2}(\Omega)}  \leq \Phi(h),
\end{equation}
where
\begin{equation}\label{phideh}
\Phi(h) \doteq \sqrt{2L} \left( \frac{1 + M}{L-a} \,  \| h \|_{L^{2}(\Gamma_{I})}
+  \| \nabla h \|_{L^{2}(\Gamma_{I})} + M \, \left\| \dfrac{\partial h_{2}}{\partial y} + \dfrac{\partial h_{3}}{\partial z} \right\|_{L^{2}(\Gamma_{I})} \right).
\end{equation}
\end{theorem}
\noindent
\begin{proof}
Consider the vector field $A_{1} \in H^{1}(Q)$ defined by
$$
A_{1}(x,y,z) = h(y,z) \qquad \text{for a.e.} \ \ (x,y,z) \in Q,
$$
so that $A_{1} = h$ on $\Gamma_{I} \cup \Gamma_{O}$, $A_{1}$ vanishes on $([-L,L]^2 \times \{\pm L\}) \cup ([-L,L] \times \{\pm L\} \times [-L,L])$ and
\begin{equation}\label{normasa1}
\begin{aligned}
\| A_1 \|_{L^{2}(Q)} & = \sqrt{2L} \, \| h \|_{L^{2}(\Gamma_{I})} \, , \qquad \| \nabla A_1 \|_{L^{2}(Q)} = \sqrt{2L} \, \| \nabla h \|_{L^{2}(\Gamma_{I})} \, , \\[6pt]
& \| \nabla \cdot A_1 \|_{L^{2}(Q)} = \sqrt{2L} \, \left\| \dfrac{\partial h_{2}}{\partial y} + \dfrac{\partial h_{3}}{\partial z} \right\|_{L^{2}(\Gamma_{I})} \, .
\end{aligned}
\end{equation}
Moreover, since the outward unit normal to $\Gamma_{I}$ and $\Gamma_{O}$ is respectively $\mp e_{1}$, we have
$$
\int_{\partial Q} A_{1} \cdot \hat{n} = - \int_{(-L,L)^2} h(y,z) \cdot e_{1} \, dy \, dz + \int_{(-L,L)^2} h(y,z) \cdot e_{1} \, dy \, dz = 0.
$$
We are then in position to invoke \cite[Theorem 2.1]{gazspefra}, deducing the existence of a vector field $v_0\in H^{1}(\Omega)$ satisfying \eqref{div} such that
\begin{equation}\label{ilaria}
\| \nabla v_{0} \|_{L^{2}(\Omega)}  \leq \frac{1 + M}{L-a} \,  \| A_{1} \|_{L^{2}( Q)}
+  \| \nabla A_{1} \|_{L^{2}(Q)} + M \, \| \nabla \cdot A_{1} \|_{L^{2}(Q)} \, .
\end{equation}
The proof is complete after inserting \eqref{normasa1} into \eqref{ilaria}.
\end{proof}

\section{Stokes and Navier-Stokes boundary-value problems} \label{boundaryvalueproblems}

Let us start by explaining why the do-nothing boundary condition \eqref{nsstokes0}$_4$ arises naturally while studying the weak formulation associated to problem \eqref{nsstokes0}; indeed, assume that $u \in \mathcal{C}^{2}(\overline{\Omega})$ and $p \in \mathcal{C}^{1}(\overline{\Omega})$ solve \eqref{nsstokes0} in the classical sense. We take a vector function $\varphi \in H^{1}(\Omega)$ and integrate by parts the equation of conservation of momentum \eqref{nsstokes0}$_1$ in the following way:
\begin{equation} \label{testing1}
\begin{aligned}
& \int_{\Omega} \left[ -\eta\Delta u+(u\cdot\nabla)u+\nabla p \right] \cdot \varphi  \\[6pt]
& \hspace{-5mm} = \eta \int_{\Omega} \nabla u \cdot \nabla \varphi  - \eta \int_{\partial \Omega}  (\nabla u \cdot \hat{n}) \cdot \varphi  + \int_{\Omega} (u \cdot \nabla)u \cdot \varphi  + \int_{\partial \Omega}  p \hat{n} \cdot \varphi  - \int_{\Omega} p  (\nabla \cdot \varphi)  \\[6pt]
& \hspace{-5mm} = \eta \int_{\Omega} \nabla u \cdot \nabla \varphi + \int_{\Omega} (u \cdot \nabla)u \cdot \varphi - \int_{\Omega} p  (\nabla \cdot \varphi) + \int_{\partial \Omega} \left( p \hat{n} - \eta \dfrac{\partial u}{\partial n} \right) \cdot \varphi \, .
\end{aligned}
\end{equation}
If, in addition, we assume that $\varphi$ is divergence-free and vanishes on $\Gamma_{I} \cup \Gamma_{W}$, from \eqref{testing1} we obtain
$$
0 = \int_{\Omega} \left[ -\eta\Delta u+(u\cdot\nabla)u+\nabla p \right] \cdot \varphi = \eta \int_{\Omega} \nabla u \cdot \nabla \varphi + \int_{\Omega} (u \cdot \nabla)u \cdot \varphi + \int_{\Gamma_{O}} \left( p \hat{n} - \eta \dfrac{\partial u}{\partial n} \right) \cdot \varphi, 
$$
so that the boundary condition \eqref{nsstokes0}$_{4}$ yields the equality
$$
\eta \int_{\Omega} \nabla u \cdot \nabla \varphi + \int_{\Omega} (u \cdot \nabla)u \cdot \varphi  = 0 \, ,
$$
where it suffices to have $u \in H^{1}(\Omega)$ and the scalar pressure is no longer present. Notice that since $u$ vanishes on $\Gamma_{W}$, the incompressibility condition in \eqref{nsstokes0}$_{1}$ also implies the identity
$$
\int_{\Gamma_{I}}  h \cdot e_{1} = \int_{\Gamma_{O}} u \cdot e_{1} \, . 
$$
This motivates the introduction of the functional spaces (of vector fields) that will be employed hereafter:
$$
\begin{aligned}
\mathcal{V}_{*}(\Omega) & = \{ v \in V_{*}(\Omega) \ | \ \nabla \cdot v=0 \ \ \mbox{in} \ \ \Omega \}  \qquad \text{(for the test functions)}, \\[3pt]
\mathcal{V}(\Omega) & = \{ v \in V(\Omega) \ | \ \nabla \cdot v=0 \ \ \mbox{in} \ \ \Omega \} \qquad \text{(for the solutions of \eqref{nsstokes0})}, 
\end{aligned}
$$
which are Hilbert spaces if endowed with the scalar product $(u,v)\mapsto(\nabla u,\nabla v)_{L^2(\Omega)}$. 

\subsection{Estimates for the Stokes operator under mixed boundary conditions}
Let us consider the following boundary-value problem associated to the Stokes equations in $\Omega$:
\begin{equation}\label{stokes0}
\left\{
\begin{aligned}
& -\eta\Delta u + \nabla p= g,\ \quad  \nabla\cdot u=0 \ \ \mbox{ in } \ \ \Omega, \\[3pt]
& u=h \ \ \mbox{ on } \ \ \Gamma_{I}, \\[3pt]
& u=0 \ \ \mbox{ on } \ \ \Gamma_{W}, \\[3pt]
& -p \hat{n} + \eta \dfrac{\partial u}{\partial n} = 0 \ \ \mbox{ on } \ \ \Gamma_{O} \, ,
\end{aligned}
\right.
\end{equation}
where, as in \eqref{nsstokes0}, $u : \Omega \longrightarrow \mathbb{R}^3$ is the velocity field, $p : \Omega \longrightarrow \mathbb{R}$ is the scalar pressure, $g : \Omega \longrightarrow \mathbb{R}^3$ denotes an external forcing term, $\eta>0$ is the (constant) kinematic viscosity and the datum $h : \Gamma_{I} \longrightarrow \mathbb{R}^3$ describes the inflow boundary condition on the wall $\Gamma_{I}$. The purpose of this section is to obtain explicit estimates for the weak solution of \eqref{stokes0} (in terms of the domain and on the data), which are defined as:
\begin{definition}\label{weaksolutionstokes}
Given $g \in L^{3/2}(\Omega)$ and $h \in H_{0}^{1}(\Gamma_{I})$, we say that a vector field $u \in \mathcal{V}(\Omega)$ is a \textbf{weak solution} of \eqref{stokes0} if $u$ verifies \eqref{stokes0}$_2$ in the trace sense and
\begin{equation} \label{nstokesdebil}
\eta \int_{\Omega} \nabla u \cdot \nabla \varphi = \int_{\Omega} g \cdot \varphi \qquad \forall \varphi \in \mathcal{V}_{*}(\Omega).
\end{equation}
\end{definition}

The well-posedness of the Stokes problem \eqref{stokes0} is established in the following result, whose proof can be found in \cite[Section 3.1]{fursikov2009optimal}:
\begin{theorem} \label{wellposstokes}
Given any $g \in L^{3/2}(\Omega)$ and $h \in H^{1}_{0}(\Gamma_{I})$, there exists a unique weak solution $u \in \mathcal{V}(\Omega)$ of problem \eqref{stokes0} together with a function $p \in L^{3/2}(\Omega)$ (unique up to an additive constant) such that the pair $(u,p)$ solves \eqref{stokes0}$_1$ in distributional sense. Moreover, if $(u,p) \in W^{2,3/2}(\Omega) \times W^{1,3/2}(\Omega)$, then $p$ is uniquely defined and this pair verifies the boundary conditions \eqref{stokes0}$_2$-\eqref{stokes0}$_3$-\eqref{stokes0}$_4$ pointwise.
\end{theorem}

With the aid of the results of Section \ref{funcine} we can also derive explicit bounds for the solutions of \eqref{stokes0}:
\begin{theorem} \label{estimatesstokes}
Assume \eqref{OmegaR}. For $g \in L^{3/2}(\Omega)$ and $h \in H^{1}_{0}(\Gamma_{I})$, let $u \in \mathcal{V}(\Omega)$ be the unique weak solution of \eqref{stokes0}. We then have
\begin{equation}\label{stokesapriori1}
 \| \nabla u \|_{L^{2}(\Omega)}  \leq \sqrt{2} \left( \Phi(h) + \dfrac{1}{\sqrt{\mathcal{S}_{3}} \, \eta} \| g \|_{L^{3/2}(\Omega)} \right) ,
\end{equation}
where $\mathcal{S}_{3} > 0$ and $\Phi(h)$ are defined as in \eqref{sobolevconstants1} and \eqref{phideh}, respectively.
\end{theorem}
\noindent
\begin{proof}
In virtue of Theorem \ref{teononconstant}, let us consider a (solenoidal) vector field $v_0\in H^{1}(\Omega)$ satisfying \eqref{div} together with inequality \eqref{f:gb}. It is then clear that $u - v_{0} \in \mathcal{V}_{*}(\Omega)$, so that, by testing with $\varphi = u - v_{0}$ in \eqref{nstokesdebil}, we obtain:
\begin{equation}\label{stokesapriori2}
\eta \| \nabla u \|^{2}_{L^{2}(\Omega)} = \eta \int_{\Omega} \nabla u \cdot \nabla v_{0} + \int_{\Omega} g \cdot (u-v_{0}) \, .
\end{equation}
We apply Young's inequality to obtain
\begin{equation}\label{stokesapriori3}
\left| \int_{\Omega} \nabla u \cdot \nabla v_{0} \right| \leq \| \nabla u \|_{L^{2}(\Omega)} \| \nabla v_{0} \|_{L^{2}(\Omega)} \leq \dfrac{1}{4} \| \nabla u \|^{2}_{L^{2}(\Omega)} + \| \nabla v_{0} \|^{2}_{L^{2}(\Omega)} \, .
\end{equation}
On the other hand, in view of the embedding $V_{*}(\Omega) \subset L^{3}(\Omega)$ and of Young's inequality we have
\begin{equation}\label{stokesapriori4}
\begin{aligned}
\left| \int_{\Omega} g \cdot (u-v_{0}) \right| & \leq \| g \|_{L^{3/2}(\Omega)} \| u - v_{0} \|_{L^{3}(\Omega)} \leq \dfrac{1}{\sqrt{\mathcal{S}_{3}}} \| g \|_{L^{3/2}(\Omega)} \| \nabla (u - v_{0}) \|_{L^{2}(\Omega)} \\[3pt]
& \leq \dfrac{\eta}{4} \| \nabla u \|^{2}_{L^{2}(\Omega)} + \dfrac{1}{\mathcal{S}_{3} \, \eta} \| g \|^{2}_{L^{3/2}(\Omega)} + \dfrac{1}{\sqrt{\mathcal{S}_{3}}} \| g \|_{L^{3/2}(\Omega)} \| \nabla v_{0} \|_{L^{2}(\Omega)} \, .
\end{aligned}
\end{equation}
By inserting \eqref{stokesapriori3} and \eqref{stokesapriori4} into \eqref{stokesapriori2} we infer that
$$
\dfrac{\eta}{2} \| \nabla u \|^{2}_{L^{2}(\Omega)} \leq \eta \| \nabla v_{0} \|^{2}_{L^{2}(\Omega)} + \dfrac{1}{\mathcal{S}_{3} \, \eta} \| g \|^{2}_{L^{3/2}(\Omega)} + \dfrac{1}{\sqrt{\mathcal{S}_{3}}} \| g \|_{L^{3/2}(\Omega)} \| \nabla v_{0} \|_{L^{2}(\Omega)} \, ,
$$
which, in view of \eqref{f:gb}, yields
$$
\dfrac{1}{2} \| \nabla u \|^{2}_{L^{2}(\Omega)} \leq \Phi(h)^{2} + \dfrac{1}{\mathcal{S}_{3} \, \eta^{2}} \| g \|^{2}_{L^{3/2}(\Omega)} + \dfrac{\Phi(h)}{\sqrt{\mathcal{S}_{3}} \, \eta} \| g \|_{L^{3/2}(\Omega)} \leq \left( \Phi(h) + \dfrac{1}{\sqrt{\mathcal{S}_{3}} \, \eta} \| g \|_{L^{3/2}(\Omega)} \right)^{2} \, .
$$
\end{proof}

We can now give the following definition for the weak solutions of problem \eqref{nsstokes0}:
\begin{definition}\label{weaksolution}
Given $h \in H_{0}^{1}(\Gamma_{I})$, we say that a vector field $u \in \mathcal{V}(\Omega)$ is a \textbf{weak solution} of \eqref{nsstokes0} if $u$ verifies \eqref{nsstokes0}$_2$ in the trace sense and
\begin{equation} \label{nstokesdebil1}
\eta \int_{\Omega} \nabla u \cdot \nabla \varphi + \int_{\Omega} (u \cdot \nabla)u \cdot \varphi = 0 \qquad \forall \varphi \in \mathcal{V}_{*}(\Omega).
\end{equation}
\end{definition}

We are now in position to prove one of the main results of the present article:

\begin{theorem} \label{mainresult}
Assume \eqref{OmegaR}. Define $\mathcal{S}_{3}, \mathcal{S}_{6}, \mathcal{J}_{6} > 0$ as in \eqref{sobolevconstants1}, and let $h \in H^{1}_{0}(\Gamma_{I})$ be such that
\begin{equation}\label{threshold}
\Phi(h) \leq \dfrac{\eta^{2}}{2 \sqrt{2}} \, \dfrac{\sqrt{\mathcal{J}_{6}} \, \mathcal{S}_{3}}{(1 + \sqrt{\mathcal{S}_{3}} \, \eta)^{2}} \, ,
\end{equation} 
with $\Phi(h)$ given as in \eqref{phideh}. Then, there exists a unique weak solution $u \in \mathcal{V}(\Omega)$ of \eqref{nsstokes0} such that
\begin{equation}\label{main1}
\| \nabla u \|_{L^{2}(\Omega)} \leq \sqrt{2} \left[ 1 + \dfrac{2}{\sqrt{\mathcal{J}_{6} \, \mathcal{S}_{3}} \, \eta} \, \dfrac{\sqrt{2} \, (\sqrt{\mathcal{S}_{6}} + \sqrt{\mathcal{J}_{6}}) \, (1 + \sqrt{\mathcal{S}_{3}} \, \eta) \, \Phi(h)^{2}}{ \sqrt{\mathcal{J}_{6} \mathcal{S}_{6}} \, \mathcal{S}_{3} \, \eta^{2} - \sqrt{2} \, (\sqrt{\mathcal{S}_{6}} + \sqrt{\mathcal{J}_{6}}) \, (1 + \sqrt{\mathcal{S}_{3}} \, \eta) \, \Phi(h)} \right] \Phi(h) \, .
\end{equation}
\end{theorem}
\noindent
\begin{proof}
Following \cite[Theorem 3.7]{fursikov2009optimal}, for any $p \in (1,\infty)$ we set
$$
L_{*}^{p}(\Omega) = \{ v \in L^{p}(\Omega) \ | \ \nabla \cdot v=0 \ \ \mbox{in} \ \ \Omega, \ \ v \cdot \hat{n} =0 \ \ \mbox{on} \ \ \Gamma_{I} \cup \Gamma_{W} \} \, ,
$$
which corresponds to the closure of $\mathcal{V}_{*}(\Omega)$ with respect to the $L^{p}(\Omega)$-norm (see also \cite{maslennikova1983approximation}), so that  $L_{*}^{p}(\Omega)$ is a closed subspace of $L^{p}(\Omega)$. Notice that for every $v \in L_{*}^{p}(\Omega)$, the normal component of its trace is well-defined as an element of $W^{-1/p, \, p}(\partial \Omega)$. The \textit{projection operator} $\mathcal{P} : L^{3/2}(\Omega) \longrightarrow L_{*}^{3/2}(\Omega)$ is introduced in such a way that, given any $g \in L^{3/2}(\Omega)$, $\mathcal{P}(g) \in L_{*}^{3/2}(\Omega)$ is the unique solution of
	\begin{equation}\label{project1}
	\int_{\Omega} \mathcal{P}(g) \cdot \varphi = \int_{\Omega} g \cdot \varphi  \qquad \forall \varphi \in L_{*}^{3}(\Omega).
	\end{equation}
We then define the operator $\widehat{\mathcal{P}} : L^{3/2}_{*}(\Omega) \longrightarrow L^{3/2}_{*}(\Omega)$ by
\begin{equation} \label{operatorpgorro}
\widehat{\mathcal{P}}(g) = - \mathcal{P}\left( (u \cdot \nabla)u \right ) \qquad \forall g \in L_{*}^{3/2}(\Omega),
\end{equation}	
where, given $g \in L^{3/2}_{*}(\Omega)$, $u \in \mathcal{V}(\Omega)$ is the unique weak solution of the Stokes problem \eqref{stokes0} with a right-hand side given by $g$. Notice that, in view of the embedding $H^{1}(\Omega) \subset L^{6}(\Omega)$, we have that $(u \cdot \nabla)u \in L^{3/2}(\Omega)$ together with
\begin{equation} \label{norm32}
\| (u \cdot \nabla) u \|_{L^{3/2}(\Omega)} \leq \| u \|_{L^{6}(\Omega)} \| \nabla u \|_{L^{2}(\Omega)} \leq \dfrac{1}{\sqrt{\mathcal{J}_{6}}} \| \nabla u \|^{2}_{L^{2}(\Omega)} \, ,
\end{equation}
see \eqref{sobolevconstants11}, so that $\widehat{\mathcal{P}}$ is correctly defined. Now, let us suppose that $\hat{\mathcal{P}}$ has a fixed point $\hat{g} \in L_{*}^{3/2}(\Omega)$: denoting by $\hat{u} \in \mathcal{V}(\Omega)$ is the unique weak solution of \eqref{stokes0} (with a right-hand side $\hat{g}$), by definition of $\mathcal{P}$ and by Definition \ref{weaksolutionstokes} we would have
$$
\eta \int_{\Omega} \nabla \hat{u} \cdot \nabla \varphi = \int_{\Omega} \hat{g} \cdot \varphi = - \int_{\Omega}  (\hat{u} \cdot \nabla)\hat{u}  \cdot \varphi \qquad \forall \varphi \in \mathcal{V}_{*}(\Omega) \subset L_{*}^{3}(\Omega),
$$
so that, in particular, $\hat{u}$ is a weak solution of \eqref{nsstokes0} as in Definition \ref{weaksolution}. Therefore, the number of fixed points of $\widehat{\mathcal{P}}$ equals the number of weak solutions of \eqref{nsstokes0}, and in the remaining part of the proof we will show that $\widehat{\mathcal{P}}$ possesses one and only one fixed point in the ball
$$
\mathcal{B}_{0} \doteq \left\{ g \in L_{*}^{3/2}(\Omega) \ \big| \ \| g \|_{L^{3/2}(\Omega)} \leq \Phi(h) \right\} \, .
$$
In view of Theorem \ref{estimatesstokes}, \eqref{norm32} and hypothesis \eqref{threshold}, for every $g \in \mathcal{B}_{0}$ we have
\begin{equation} \label{sperone4antes}
\begin{aligned}
\| \widehat{\mathcal{P}}(g) \|_{L^{3/2}(\Omega)} & =  \| \mathcal{P} ((u \cdot \nabla) u) \|_{L^{3/2}(\Omega)} = \| (u \cdot \nabla) u \|_{L^{3/2}(\Omega)} \leq \dfrac{1}{\sqrt{\mathcal{J}_{6}}} \| \nabla u \|^{2}_{L^{2}(\Omega)} \\[6pt]
& \leq \dfrac{2}{\sqrt{\mathcal{J}_{6}}} \left( \Phi(h) + \dfrac{1}{\sqrt{\mathcal{S}_{3}} \, \eta} \| g \|_{L^{3/2}(\Omega)} \right)^{2} \leq \dfrac{2}{\sqrt{\mathcal{J}_{6}}} \left( 1 + \dfrac{1}{\sqrt{\mathcal{S}_{3}} \, \eta} \right)^{2} \, \Phi(h)^2 \leq \dfrac{\Phi(h)}{\sqrt{2}} \, ,
\end{aligned}
\end{equation}
so that $\widehat{\mathcal{P}}$ maps $\mathcal{B}_{0}$ into itself. Now, given $g_{1}, g_{2} \in \mathcal{B}_{0}$, we denote by $u_{1}, u_{2} \in \mathcal{V}(\Omega)$, respectively, the unique weak solutions of \eqref{stokes0} with right-hand sides given by $g_{1}$ and $g_{2}$. We then have
\begin{equation} \label{sperone4}
\begin{aligned}
\| \widehat{\mathcal{P}}(g_{1}) - \widehat{\mathcal{P}}(g_{2}) \|_{L^{3/2}(\Omega)} & = \| \mathcal{P}\left( (u_{2} \cdot \nabla)u_{2} - (u_{1} \cdot \nabla)u_{1} \right) \|_{L^{3/2}(\Omega)} \\[6pt]
& = \| (u_{2} \cdot \nabla)u_{2} - (u_{1} \cdot \nabla)u_{1} \|_{L^{3/2}(\Omega)} \\[6pt]
& \leq \| (u_{2} \cdot \nabla) (u_{2} - u_{1}) \|_{L^{3/2}(\Omega)} + \| ((u_{2} - u_{1}) \cdot \nabla)u_{1} \|_{L^{3/2}(\Omega)} \\[6pt]
& \leq \| u_{2} \|_{L^{6}(\Omega)} \| \nabla (u_{2} - u_{1}) \|_{L^{2}(\Omega)} + \| u_{2} - u_{1} \|_{L^{6}(\Omega)} \| \nabla u_{1} \|_{L^{2}(\Omega)} \\[6pt]
& \leq \left( \dfrac{1}{\sqrt{\mathcal{J}_{6}}} \| \nabla u_{2} \|_{L^{2}(\Omega)} + \dfrac{1}{\sqrt{\mathcal{S}_{6}}} \| \nabla u_{1} \|_{L^{2}(\Omega)} \right) \| \nabla ( u_{2} - u_{1}) \|_{L^{2}(\Omega)} \, ,
\end{aligned} 
\end{equation}
see again \eqref{sobolevconstants11}, noticing that $u_{2} - u_{1} \in \mathcal{V}_{*}(\Omega)$. Therefore, we test with $\varphi = u_{2} - u_{1}$ in both the weak formulations \eqref{nstokesdebil} satisfied by $u_{1}$ and $u_{2}$, and then subtract the obtained identities to deduce
$$
\eta \| \nabla ( u_{2} - u_{1}) \|^{2}_{L^{2}(\Omega)} = (g_{2} - g_{1}, u_{2} - u_{1})_{L^{2}(\Omega)} \leq \| g_{2} - g_{1} \|_{L^{3/2}(\Omega)} \| u_{2} - u_{1} \|_{L^{3}(\Omega)} \, ,
$$
so that, in view of \eqref{sobolevconstants11},
\begin{equation} \label{sperone5}
\| \nabla ( u_{2} - u_{1}) \|_{L^{2}(\Omega)} \leq \dfrac{1}{\sqrt{\mathcal{S}_{3}} \, \eta} \| g_{2} - g_{1} \|_{L^{3/2}(\Omega)} \, .
\end{equation}
On the other hand, from Theorem \ref{estimatesstokes} we have the estimates
\begin{equation}\label{sperone6}
\| \nabla u_{i} \|_{L^{2}(\Omega)}  \leq \sqrt{2} \left( \Phi(h) + \dfrac{1}{\sqrt{\mathcal{S}_{3}} \, \eta} \| g_{i} \|_{L^{3/2}(\Omega)} \right) \leq \sqrt{2} \left( 1 + \dfrac{1}{\sqrt{\mathcal{S}_{3}} \, \eta} \right) \Phi(h), \qquad \text{for} \ \ i \in \{1,2\} .
\end{equation}
Because of \eqref{sperone5}-\eqref{sperone6}, inequality \eqref{sperone4} becomes
$$
\| \widehat{\mathcal{P}}(g_{1}) - \widehat{\mathcal{P}}(g_{2}) \|_{L^{3/2}(\Omega)} \leq \dfrac{\sqrt{2}}{\sqrt{\mathcal{S}_{3}} \, \eta}  \left( \dfrac{1}{\sqrt{\mathcal{J}_{6}}} + \dfrac{1}{\sqrt{\mathcal{S}_{6}}} \right) \left( 1 + \dfrac{1}{\sqrt{\mathcal{S}_{3}} \, \eta} \right) \Phi(h) \, \| g_{1} - g_{2} \|_{L^{3/2}(\Omega)} \qquad \forall g_{1}, g_{2} \in \mathcal{B}_{0} \, ,
$$
where, owing to \eqref{threshold} and since $\mathcal{J}_{6} \leq \mathcal{S}_{6}$,
\begin{equation} \label{beta}
\beta \doteq \dfrac{\sqrt{2}}{\sqrt{\mathcal{S}_{3}} \, \eta}  \left( \dfrac{1}{\sqrt{\mathcal{J}_{6}}} + \dfrac{1}{\sqrt{\mathcal{S}_{6}}} \right) \left( 1 + \dfrac{1}{\sqrt{\mathcal{S}_{3}} \, \eta} \right) \Phi(h) < \dfrac{2 \sqrt{2}}{\sqrt{\mathcal{J}_{6}}} \left( 1 + \dfrac{1}{\sqrt{\mathcal{S}_{3}} \, \eta} \right)^{2} \, \Phi(h) \leq 1 \, .
\end{equation}
The operator $\widehat{\mathcal{P}}$ then satisfies the contraction property on $\mathcal{B}_{0}$ and, therefore, admits a unique fixed point $\hat{g} \in \mathcal{B}_{0}$ which is the limit of the sequence $\{ g_{n} \}_{n \in \mathbb{N}} \subset \mathcal{B}_{0}$ defined inductively as
$$
g_{0} \doteq 0, \qquad g_{n} = \widehat{\mathcal{P}}(g_{n-1}) = \sum_{k=1}^{n} (g_{k} - g_{k-1})  \qquad \forall n \geq 1.
$$
As a consequence, we can estimate the norm of $\hat{g}$ in the following way
\begin{equation} \label{sperone8}
\begin{aligned}
\| \hat{g} \|_{L^{3/2}(\Omega)} & \leq \lim_{n \to \infty} \sum_{k=1}^{n} \| g_{k} - g_{k-1} \|_{L^{3/2}(\Omega)} \leq \sum_{k=0}^{\infty} \beta^{k} \, \| \widehat{\mathcal{P}}(0) \|_{L^{3/2}(\Omega)} = \dfrac{\beta}{1 - \beta} \, \| \mathcal{P}\left( (u_{0} \cdot \nabla)u_{0} \right) \|_{L^{3/2}(\Omega)} \\[6pt]
& = \dfrac{\beta}{1 - \beta} \, \| (u_{0} \cdot \nabla)u_{0} \|_{L^{3/2}(\Omega)} \leq \dfrac{2}{\sqrt{\mathcal{J}_{6}}} \, \dfrac{\beta}{1 - \beta} \, \Phi(h)^2 \, ,
\end{aligned}
\end{equation}
with $\beta \in (0,1)$ as in \eqref{beta}, $u_{0} \in \mathcal{V}(\Omega)$ the unique weak solution of \eqref{stokes0} with $g=0$ and the last inequality in \eqref{sperone8} follows from combining \eqref{norm32} with \eqref{stokesapriori1}. By construction, the function $\hat{u} \in \mathcal{V}$ (i.e., the unique weak solution of \eqref{stokes0} with right-hand side given by $\hat{g}$) is the unique weak solution of \eqref{nsstokes0}. Theorem \ref{estimatesstokes} and \eqref{sperone8} then yield the following estimate:
\begin{equation} \label{estimarsol}
\| \nabla \hat{u} \|_{L^{2}(\Omega)} \leq \sqrt{2} \left( \Phi(h) + \dfrac{1}{\sqrt{\mathcal{S}_{3}} \, \eta} \| \hat{g} \|_{L^{3/2}(\Omega)} \right) \leq \sqrt{2} \left( 1 + \dfrac{2}{\sqrt{\mathcal{J}_{6} \, \mathcal{S}_{3}} \, \eta} \, \dfrac{\beta}{1 - \beta} \, \Phi(h) \right) \Phi(h) \, ,
\end{equation}
which concludes the proof.
\end{proof}

\begin{remark} \label{remark1}
It can be proved (see \cite[Theorem 4.1]{fursikov2009optimal}) that, under hypothesis \eqref{threshold}, there also exists an associated pressure $p \in L^{2}(\Omega)$ such that $(u,p) \in \mathcal{V}(\Omega) \times L^2(\Omega)$ solves \eqref{nsstokes0}$_1$ in distributional sense.
\end{remark}

\begin{remark}
	We point out that Theorems \ref{wellposstokes} and \ref{mainresult} remain valid if, instead of $H^{1}_{0}(\Gamma_{I})$, we take a boundary datum $h$ belonging to
	$$
	H^{1/2}_{00}(\Gamma_{I}) \doteq \{ h \in L^{2}(\Gamma_{I}) \ | \ (\exists \hat{h} \in H^{1}(\Omega)) \ \ \hat{h} = h \ \ \mbox{on} \ \ \Gamma_{I}, \ \ \hat{h} = 0 \ \ \mbox{on} \ \ \partial\Omega \setminus \Gamma_{I} \},
	$$
	which is a closed subspace of $H^{1/2}(\Gamma_{I})$, see \cite[Chapter VII]{dautray1999mathematical} for more details. Nevertheless, for the purpose of obtaining explicit bounds, we restrict ourselves to inflows in $H^{1}_{0}(\Gamma_{I})$, see Theorem \ref{teononconstant}, and emphasize that the estimates given in \eqref{stokesapriori1} and \eqref{main1} depend on the solenoidal extension employed. 
\end{remark}

\noindent
Combining Theorem \ref{wellposstokes} with the results of Section \ref{funcine}, we formulate the main result of the present article:

\begin{corollary} \label{maincor}
Assume \eqref{OmegaR}. If the boundary datum $h \in H^{1}_{0}(\Gamma_{I})$ is such that
\begin{equation}\label{thresholdcor}
\Phi(h) \leq \dfrac{\eta^{2}}{4 \sqrt{2} \, \pi^{1/3}} \, \dfrac{\left[ \pi^2 - \dfrac{2}{L^2} \min \left\{ \sqrt[3]{\dfrac{3}{2\pi} (|Q|-|K|)}  \, , \,  \dfrac{4L}{3} \right\}^{2} \right]^{1/2}}{\left[ \eta + \dfrac{1}{\sqrt{5}} \left( \dfrac{7^{7/3}}{3} \right)^{1/4} \left( \dfrac{1}{2 \pi^{5}} \right)^{1/6} \, \min \left\{ \sqrt[3]{\dfrac{3}{2\pi} (|Q|-|K|)}  \, , \,  \dfrac{4L}{3} \right\}^{1/2} \right]^{2}} \, ,
\end{equation} 
with $\Phi(h)$ given as in \eqref{phideh}, then problem \eqref{nsstokes0} admits a unique weak solution $u \in \mathcal{V}(\Omega)$ such that
\begin{equation}\label{maincor1}
\| \nabla u \|_{L^{2}(\Omega)} \leq \left[ \sqrt{2} + \dfrac{1}{\sqrt{5} \, \eta} \left( \dfrac{7^{7/3}}{3} \right)^{1/4} \left( \dfrac{1}{2 \pi^{5}} \right)^{1/6} \, \min \left\{ \sqrt[3]{\frac{3}{2\pi} (|Q|-|K|)}  \, , \,  \dfrac{4L}{3} \right\}^{1/2} \, \right] \Phi(h) \, .
\end{equation}
\end{corollary}
\noindent
\begin{proof}
In virtue of Theorem \ref{boundsobolev} we have, on one hand,
$$
\sqrt{\mathcal{J}_{6}} \geq \dfrac{1}{2 \pi^{1/3}} \, \left[ \pi^2 - \dfrac{2}{L^2} \min \left\{ \sqrt[3]{\frac{3}{2\pi} (|Q|-|K|)}  \, , \,  \dfrac{4L}{3} \right\}^{2} \right]^{1/2} \, .
$$
On the other hand, Corollary \ref{boundsobolev3} gives the inequality
\begin{equation}\label{sobolev3}
\dfrac{1}{\sqrt{\mathcal{S}_{3}}} \leq \dfrac{1}{\sqrt{5}} \left( \dfrac{7^{7/3}}{3} \right)^{1/4} \left( \dfrac{1}{2 \pi^{5}} \right)^{1/6} \, \min \left\{ \sqrt[3]{\frac{3}{2\pi} (|Q|-|K|)}  \, , \,  \dfrac{4L}{3} \right\}^{1/2} \, .
\end{equation}
Thus, if $h \in H^{1}_{0}(\Gamma_{I})$ is such that \eqref{thresholdcor} holds, then condition \eqref{threshold} in Theorem \ref{mainresult} is certainly satisfied, ensuring the existence of a unique solution to problem  \eqref{nsstokes0}. In the proof of Theorem \ref{mainresult}, see \eqref{sperone4antes}, we notice that the fixed point $\hat{g} \in L_{*}^{3/2}(\Omega)$ belongs to the ball of radius $\Phi(h)/\sqrt{2}$ (centered at the origin of $L_{*}^{3/2}(\Omega)$), so that from \eqref{estimarsol} we can also derive the (rougher) bound 
$$
\| \nabla u \|_{L^{2}(\Omega)} \leq \left( \sqrt{2} + \dfrac{1}{\sqrt{\mathcal{S}_{3}} \, \eta} \right) \Phi(h) \, ,
$$
which concludes the demonstration in virtue of \eqref{sobolev3}.
\end{proof}

\begin{remark}
Theorem \ref{mainresult} gives an upper bound on the ``size'' of the inlet velocity $h$ (expressed through the quantity $\Phi(h)$) that guarantees both existence and uniqueness of solutions for \eqref{nsstokes0} in some ball whose radius is given by the right-hand side of \eqref{main1}. Corollary \ref{maincor} expresses this upper bound in terms of the geometric parameters of the obstacle and domain considered. Both results should be compared with \cite[Theorem 2.1]{filippoclara} and \cite[Theorem 3.2]{gazspe}. Related topics on uniqueness and bifurcation branches for planar steady Navier-Stokes equations under Navier boundary conditions were discussed in \cite{kochuniqueness}.
\end{remark}

\section{On the fluid forces applied over the obstacle} \label{forces}

Let $(u,p)\in \mathcal{V}(\Omega) \times L^2(\Omega)$ be a
weak solution of \eqref{nsstokes0}. It is well known (see \cite[Chapter 2]{landau}) that the stress tensor of any fluid governed by \eqref{nsstokes0}$_1$ is given by the matrix
$$
\mathbb{T}(u,p) =  -p\mathbb{I}_{3} + 2 \eta e(u)  \quad \text{in} \quad \Omega, \qquad \text{with} \quad e(u) \doteq \dfrac{1}{2}[\nabla u + (\nabla u)^{\intercal}],
$$
where $\mathbb{I}_{3}$ is the $3 \times 3$-identity matrix. Then, the total force exerted by the fluid over the obstacle $K$ is formally given by
\begin{equation} \label{forcek0}
F_{K}(u,p)=- \int\limits_{\partial K} \mathbb{T}(u,p) \cdot \hat{n} \, ,
\end{equation}
where the minus sign is due to the fact that the outward unit normal $\hat{n}$ to $\Omega$ is directed towards the interior of $K$. Notice that
\eqref{forcek0} is correctly defined: since $u \in L^{6}(\Omega)$, we have that
$$
\mathbb{T}(u,p)\in L^2(\Omega)\subset L^{3/2}(\Omega) \qquad \mbox{and} \qquad \nabla\cdot \mathbb{T}(u,p)=(u \cdot \nabla)u \in L^{3/2}(\Omega).
$$
Therefore, the normal component of the trace of $\mathbb{T}(u,p)$ exists as an element of $W^{-\frac{2}{3}, \frac{3}{2}}(\partial \Omega)$, the dual space of $W^{\frac{2}{3},3}(\partial \Omega)$. Then, we can rigorously define the force as follows.

\begin{definition}\label{weakforce}
	Let $(u,p) \in \mathcal{V}(\Omega) \times L^{2}(\Omega)$ be a weak solution of \eqref{nsstokes0}. Then, the \textbf{total force} exerted by the
	fluid over the obstacle $K$ is given by
	$$
	F_{K}(u,p)= - \langle \mathbb{T}(u,p) \cdot \hat{n} , 1 \rangle_{\partial K},
	$$
	where $\langle \cdot , \cdot \rangle_{\partial K}$ denotes the duality pairing between $W^{-\frac{2}{3}, \frac{3}{2}}(\partial K)$ and $W^{\frac{2}{3}, 3}(\partial K)$. Whenever the inflow velocity $h \in H^{1}_{0}(\Gamma_{I})$ is parallel to the vector $e_{1}$, we define the \textbf{drag} and \textbf{lift} over $K$ respectively as
	$$
	\mathcal{D}_{K}(u,p) = F_{K}(u,p) \cdot e_{1} \qquad \text{and} \qquad \mathcal{L}_{K}(u,p) = F_{K}(u,p) \cdot e_{3} \, .
	$$
\end{definition}

The main goal of this section is to provide explicit upper bounds on the drag and lift exerted over $K$, and to illustrate our results through a concrete example. Firstly we derive a more handy formula for the drag and lift, in the spirit of \cite[Lemma 3]{bonheure2020equilibrium}. 

\begin{proposition} \label{newformula}
Let $(u,p) \in \mathcal{V}(\Omega) \times L^{2}(\Omega)$ be a weak solution of \eqref{nsstokes0}, and define $\mathcal{J}_{6} > 0$ as in \eqref{sobolevconstants1}. For $i \in \{1,3\}$, let $q_{i} \in H^{1}(\Omega)$ be a vector field such that
\begin{equation} \label{div2}
\nabla\cdot q_{i} = 0 \ \text{ in } \ \Omega, \qquad q_{i} = e_{i} \ \text{ on } \ \partial K, \qquad q_{i} = 0 \ \text{ on } \ \partial Q \, .
\end{equation}
We then have
\begin{equation} \label{newformula1}
\mathcal{D}_{K}(u,p) = -\int_{\Omega} \left[ 2 \eta \, e(u) \cdot \nabla q_{1} + (u \cdot \nabla) u \cdot q_{1} \right] \, , \quad \mathcal{L}_{K}(u,p) = -\int_{\Omega} \left[ 2 \eta \, e(u) \cdot \nabla q_{3} + (u \cdot \nabla) u \cdot q_{3} \right] \, ,
\end{equation}
so that
\begin{equation} \label{newformula2}
\begin{aligned}
& | \mathcal{D}_{K}(u,p) | \leq \left( 2 \eta \, \| \nabla q_{1} \|_{L^{2}(\Omega)} + \dfrac{1}{\sqrt{\mathcal{J}_{6}}} \| \nabla u \|_{L^{2}(\Omega)} \| q_{1} \|_{L^{3}(\Omega)} \right) \| \nabla u \|_{L^{2}(\Omega)} \, , \\[6pt] 
& | \mathcal{L}_{K}(u,p) | \leq \left( 2 \eta \, \| \nabla q_{3} \|_{L^{2}(\Omega)} + \dfrac{1}{\sqrt{\mathcal{J}_{6}}} \| \nabla u \|_{L^{2}(\Omega)} \| q_{3} \|_{L^{3}(\Omega)} \right) \| \nabla u \|_{L^{2}(\Omega)} \, .
\end{aligned}
\end{equation}
\end{proposition}
\noindent
\begin{proof}
First of all, for $i \in \{1,3\}$,  the existence of $q_{i} \in H^{1}(\Omega)$ satisfying \eqref{div2} follows from \cite{ladyzhenskaya1978some} (see also \cite[Proposition 2.1]{gazspe}). We now multiply \eqref{nsstokes0}$_1$ by $q_i$ and integrate by parts as in \eqref{testing1} to obtain
$$
0 = \int_{\Omega} \left[ -\eta\Delta u+(u\cdot\nabla)u+\nabla p \right] \cdot q_{i}  = \eta \int_{\Omega} \nabla u \cdot \nabla q_{i} + \int_{\Omega} (u \cdot \nabla)u \cdot q_{i} + \int_{\partial \Omega} \left( p \hat{n} - \eta \dfrac{\partial u}{\partial n} \right) \cdot q_{i} \, ,
$$
so that
\begin{equation} \label{newformula3}
e_{i} \cdot \int_{\partial K} \left( - p \mathbb{I}_{3} + \eta \, \nabla u \right) \cdot \hat{n} =  \eta \int_{\Omega} \nabla u \cdot \nabla q_{i} + \int_{\Omega} (u \cdot \nabla)u \cdot q_{i} \, .
\end{equation}
On the other hand, after integrating by parts we get
$$
\int_{\Omega} (\nabla u)^{\intercal} \cdot \nabla q_{i} =  e_{i} \cdot \int_{\partial K}  (\nabla u)^{\intercal} \cdot \hat{n} \, ,
$$
which, after being added to \eqref{newformula3}, yields both identities in \eqref{newformula1}. The estimates in \eqref{newformula2} follow from H\"older's inequality and from \eqref{norm32}:
$$
\begin{aligned}
| \mathcal{D}_{K}(u,p) | & \leq \eta \, \| 2 e(u) \|_{L^{2}(\Omega)} \| \nabla q_{1} \|_{L^{2}(\Omega)} +  \| (u \cdot \nabla) u \|_{L^{3/2}(\Omega)} \| q_{1} \|_{L^{3}(\Omega)} \\[6pt]
& \leq 2 \eta \, \| \nabla u \|_{L^{2}(\Omega)} \| \nabla q_{1} \|_{L^{2}(\Omega)} + \dfrac{1}{\sqrt{\mathcal{J}_{6}}} \| \nabla u \|^{2}_{L^{2}(\Omega)} \| q_{1} \|_{L^{3}(\Omega)} \, ,
\end{aligned}
$$
and then similarly for $\mathcal{L}_{K}(u,p)$.
\end{proof}

It is clear that the result of Proposition \ref{newformula} depends on the vector fields $q_{i} \in H^{1}(\Omega)$ that act as solenoidal extensions of the unit vectors $e_{i}$, for $i \in \{1,3\}$. Since in this case we are dealing with constant boundary data, we can build such solenoidal extensions directly with the use of cut-off functions (without employing Bogovskii's formula, as in Theorem \ref{teononconstant}), following the classical method of Ladyzhenskaya \cite[Chapter 5]{ladyzhenskaya1969mathematical}, see also \cite[Section 3]{galdi2007further}.

\begin{proposition} \label{solext}
Assume \eqref{OmegaR}. There exist two vector fields $q_{1}, q_{3} \in H^{1}(\Omega)$ satisfying \eqref{div2} together with the estimates
$$
\begin{aligned}
& \| q_{1} \|_{L^{3}(\Omega)} \leq \dfrac{\sqrt[3]{abc}}{a} \left( 1 + \dfrac{3}{4} \, \sqrt{b^2 + c^2} \, \sqrt{\dfrac{1}{a^2} + \dfrac{1}{b^2} + \dfrac{1}{c^2}} \right) (L+a) \, , \\[6pt]
& \| q_{3} \|_{L^{3}(\Omega)} \leq \dfrac{\sqrt[3]{abc}}{a} \left( 1 + \dfrac{3}{4} \, \sqrt{a^2 + b^2} \, \sqrt{\dfrac{1}{a^2} + \dfrac{1}{b^2} + \dfrac{1}{c^2}} \right) (L+a) \, ,
\end{aligned}
$$
and 
$$
\begin{aligned}
\| \nabla q_{1} \|_{L^{2}(\Omega)} & \leq 3 \sqrt{bc(L-a)} \Bigg[ \sqrt{\dfrac{1}{a^2} + \dfrac{1}{b^2} + \dfrac{1}{c^2}} + \dfrac{1}{a} + \dfrac{1}{2b} + \dfrac{1}{2c} \\[6pt]
& \hspace{4mm} + \dfrac{3}{4} \, \dfrac{L+a}{L-a} \left( \dfrac{b+c}{a} \sqrt{\dfrac{64}{9a^2} + \dfrac{1}{b^2} + \dfrac{1}{c^2}} + \sqrt{\dfrac{1}{a^2} + \dfrac{64}{9b^2} + \dfrac{1}{c^2}} + \sqrt{\dfrac{1}{a^2} + \dfrac{1}{b^2} + \dfrac{64}{9c^2}} \right) \Bigg] , \\[6pt]
\| \nabla q_{3} \|_{L^{2}(\Omega)} & \leq 3 \sqrt{bc(L-a)} \Bigg[ \sqrt{\dfrac{1}{a^2} + \dfrac{1}{b^2} + \dfrac{1}{c^2}} + \dfrac{1}{c} + \dfrac{1}{2a} + \dfrac{1}{2b} \\[6pt]
& \hspace{4mm} + \dfrac{3}{4} \, \dfrac{L+a}{L-a} \left( \dfrac{a+b}{c} \sqrt{\dfrac{1}{a^2} + \dfrac{1}{b^2} + \dfrac{64}{9c^2}} + \sqrt{\dfrac{64}{9a^2} + \dfrac{1}{b^2} + \dfrac{1}{c^2}} + \sqrt{\dfrac{1}{a^2} + \dfrac{64}{9b^2} + \dfrac{1}{c^2}} \right) \Bigg] .
\end{aligned}
$$
\end{proposition}
\noindent
\begin{proof}
For any sufficiently small $\varepsilon > 0$ we define the function $\phi_{\varepsilon} : \mathbb{R} \longrightarrow \mathbb{R}$ as
$$
\phi_{\varepsilon} (t) = \begin {cases}
0  & \text{ if } \ \ t \in (-\infty,-1-\varepsilon] \cup [1+\varepsilon,\infty)
\\ \noalign{\medskip}
\dfrac{1}{\varepsilon^{3}} \left[ 2 |t|^3 - 3(\varepsilon+2)t^{2} + 6(1+\varepsilon)|t| + \varepsilon^{3} - 3\varepsilon - 2 \right]  & \text { if } \ \ t \in (-1-\varepsilon,-1) \cup (1,1+\varepsilon)
\\ \noalign{\medskip}
1 & \text{ if } \ \ t \in [-1,1] \, ,
\end{cases}
$$
see also \cite[Section 3]{filippoclara}, whose plot for $\varepsilon = 1/2$ is displayed below:
\begin{figure}[H]
	\begin{center}
		\includegraphics[height=40mm,width=80mm]{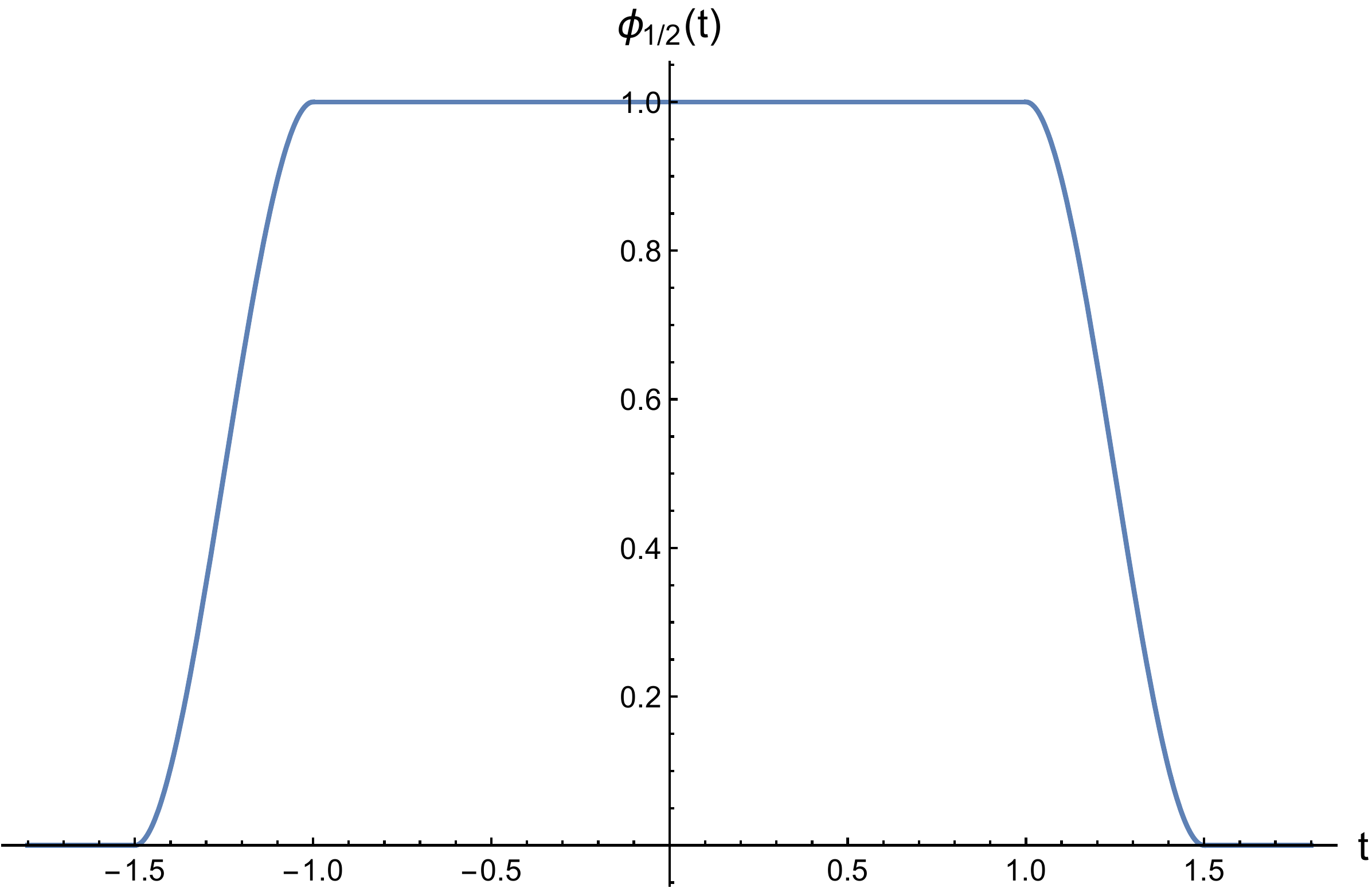} \quad \includegraphics[height=40mm,width=80mm]{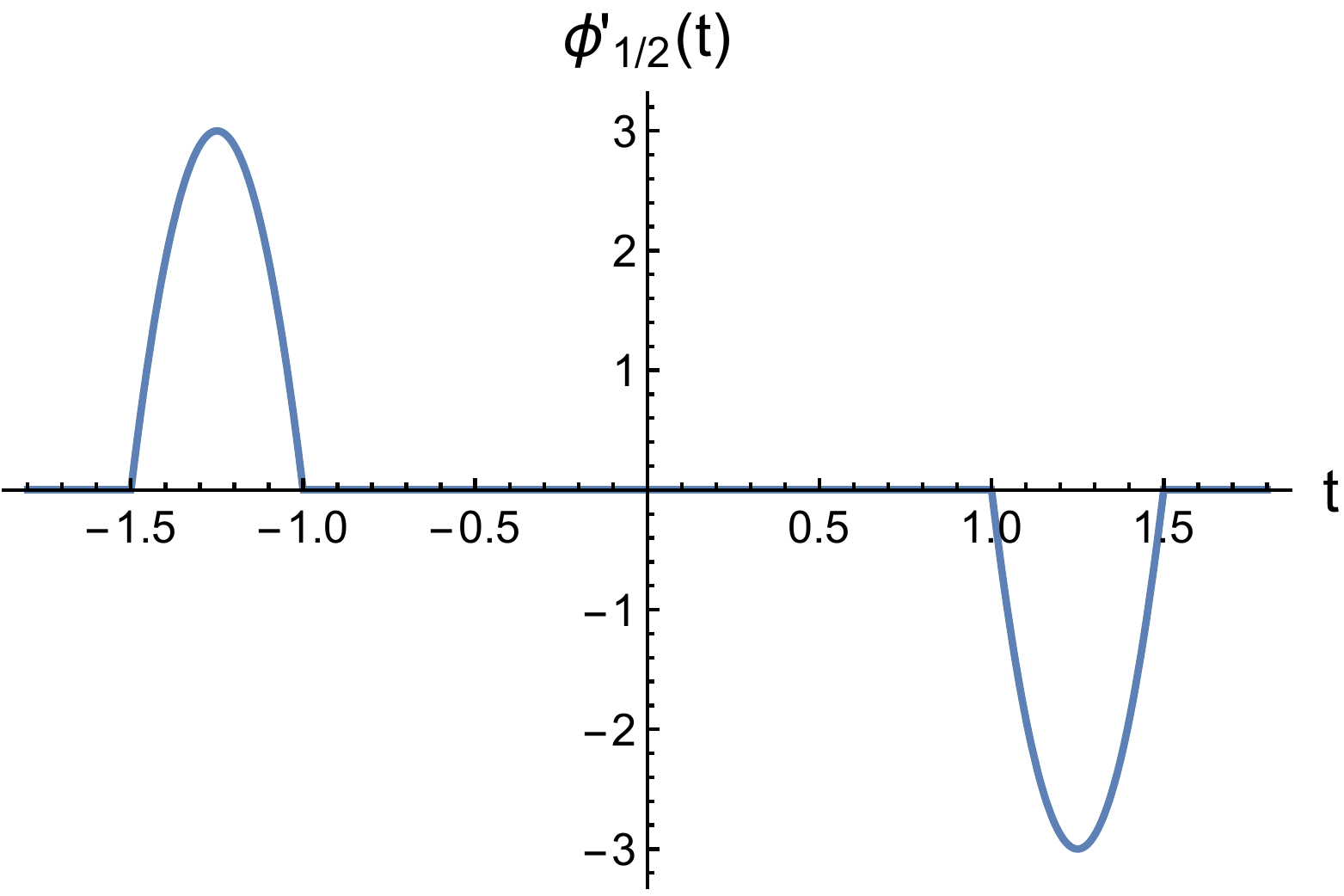}
	\end{center}
	\vspace*{-3mm}
	\caption{Graph of the function $\phi_{\varepsilon}$ (left) and of its derivative (right), for $\varepsilon = 1/2$.}\label{sym3}
\end{figure}
\noindent
Then, $\phi_{\varepsilon} \in \mathcal{C}^{1}(\mathbb{R})$, $\text{supp}(\phi_{\varepsilon}) \subset (-1-\varepsilon,1+\varepsilon)$, $\phi_{\varepsilon}'(1)=\phi_{\varepsilon}'(-1)=0$, so that  $\phi_{\varepsilon} \in H^{2}(\mathbb{R})$. In fact, we have that $\phi_{\varepsilon} \in W^{2,\infty}(\mathbb{R})$ with
\begin{equation} \label{estimar1}
\| \phi_{\varepsilon} \|_{L^{\infty}(\mathbb{R})} =1 \, , \qquad
\| \phi_{\varepsilon}' \|_{L^{\infty}(\mathbb{R})} = \dfrac{3}{2\varepsilon} \, , \qquad \| \phi_{\varepsilon}'' \|_{L^{\infty}(\mathbb{R})} = \dfrac{6}{\varepsilon^{2}} \, .
\end{equation}
We  now chose any $\varepsilon > 0$ such that, say,
$$
\varepsilon \leq \dfrac{1}{2} \min \left\{ \dfrac{L-a}{a} \, , \dfrac{L-b}{b} \, , \dfrac{L-c}{c} \, \right\} = \dfrac{L-a}{2a} \, ,
$$
and define the function $\xi_{\varepsilon} : \mathbb{R}^3 \longrightarrow \mathbb{R}$ as
$$
\xi_{\varepsilon}(x,y,z) = \phi_{\varepsilon}\left( \dfrac{x}{a} \right) \phi_{\varepsilon}\left( \dfrac{y}{b} \right) \phi_{\varepsilon}\left( \dfrac{z}{c} \right)  \qquad \forall (x,y,z) \in \mathbb{R}^3 \, ,
$$
so that $\xi_{\varepsilon} \in \mathcal{C}^{1}(\mathbb{R}^3)$, $\xi_{\varepsilon} \equiv 1$ in $\overline{P}$ and
$$
\text{supp}(\xi_{\varepsilon}) \subset (-(1+\varepsilon)a,(1+\varepsilon)a) \times (-(1+\varepsilon)b,(1+\varepsilon)b) \times (-(1+\varepsilon)c,(1+\varepsilon)c) \subset Q \setminus \overline{P} \, .
$$
The properties of $\phi_{\varepsilon}$ and \eqref{estimar1} imply that $\xi_{\varepsilon} \in H^{2}(\Omega)$ together with
\begin{equation}\label{estimar2}
\begin{aligned}
& \| \xi_{\varepsilon} \|_{L^{\infty}(\Omega)} =1 \, , \qquad
\| \nabla \xi_{\varepsilon} \|_{L^{\infty}(\Omega)} \leq \dfrac{3}{2\varepsilon} \sqrt{\dfrac{1}{a^2} + \dfrac{1}{b^2} + \dfrac{1}{c^2}}  \, , \qquad \left\|  \dfrac{\partial \xi_{\varepsilon}}{\partial x} \right\|_{L^{\infty}(\Omega)} \leq \dfrac{3}{2a\varepsilon} \, , \\[6pt]
& \left\|  \dfrac{\partial \xi_{\varepsilon}}{\partial y} \right\|_{L^{\infty}(\Omega)} \leq \dfrac{3}{2b\varepsilon} \, , \qquad \left\|  \dfrac{\partial \xi_{\varepsilon}}{\partial z} \right\|_{L^{\infty}(\Omega)} \leq \dfrac{3}{2c\varepsilon} \, , \qquad
\left\| \nabla \left( \dfrac{\partial \xi_{\varepsilon}}{\partial x} \right) \right\|_{L^{\infty}(\Omega)} \leq \dfrac{9}{4a\varepsilon^{2}} \sqrt{\dfrac{64}{9a^2} + \dfrac{1}{b^2} + \dfrac{1}{c^2}} \, , \\[6pt]
& \left\| \nabla \left( \dfrac{\partial \xi_{\varepsilon}}{\partial y} \right) \right\|_{L^{\infty}(\Omega)} \leq \dfrac{9}{4b\varepsilon^{2}} \sqrt{\dfrac{1}{a^2} + \dfrac{64}{9b^2} + \dfrac{1}{c^2}} \, , \qquad \left\| \nabla \left( \dfrac{\partial \xi_{\varepsilon}}{\partial z} \right) \right\|_{L^{\infty}(\Omega)} \leq \dfrac{9}{4c\varepsilon^{2}} \sqrt{\dfrac{1}{a^2} + \dfrac{1}{b^2} + \dfrac{64}{9c^2}} \, .
\end{aligned}
\end{equation}A straightforward computation shows that the vector fields
$$
q_{i}(x,y,z) \doteq \dfrac{1}{2} \nabla \times (\xi_{\varepsilon}(x,y,z) (e_{i} \times (x,y,z))) \qquad \forall (x,y,z) \in \overline{\Omega}\, , \ \ i \in \{1,3\}\, ,
$$
belong to $H^{1}(\Omega)$ and satisfy \eqref{div2}. More precisely, we have
$$
\left\{
\begin{aligned}
& q_{1}(x,y,z) = \xi_{\varepsilon}(x,y,z) \, e_{1} + \dfrac{1}{2} \nabla \xi_{\varepsilon}(x,y,z) \times (0,-z,y) \qquad \forall (x,y,z) \in \overline{\Omega}\, , \\[6pt]
& q_{3}(x,y,z) = \xi_{\varepsilon}(x,y,z) \, e_{3} + \dfrac{1}{2} \nabla \xi_{\varepsilon}(x,y,z) \times (-y,x,0) \qquad \forall (x,y,z) \in \overline{\Omega}\, ,
\end{aligned}
\right.
$$
and thus, by recalling that $\nabla \xi_{\varepsilon} \equiv 0$ in $\overline{P}$, we obtain the estimate
\begin{equation} \label{norml3q1}
\begin{aligned}
\| q_{1} \|_{L^{3}(\Omega)} & \leq \| \xi_{\varepsilon} \|_{L^{3}(\Omega)} + \dfrac{1+\varepsilon}{2} \sqrt{b^2 + c^2} \, \| \nabla \xi_{\varepsilon} \|_{L^{3}(\Omega)} \\[6pt]
& \leq \sqrt[3]{8 abc (1 + \varepsilon)^{3}} \, \| \xi_{\varepsilon} \|_{L^{\infty}(\Omega)} +  \dfrac{1+\varepsilon}{2} \sqrt{b^2 + c^2} \, \sqrt[3]{8 abc \varepsilon^{3}} \, \| \nabla \xi_{\varepsilon} \|_{L^{\infty}(\Omega)} \, ,
\end{aligned}
\end{equation}
and in a similar way:
\begin{equation} \label{norml3q3}
\| q_{3} \|_{L^{3}(\Omega)} \leq 2 \, \sqrt[3]{abc} \, (1+\varepsilon) \left( \| \xi_{\varepsilon} \|_{L^{\infty}(\Omega)} +  \dfrac{\varepsilon}{2} \sqrt{a^2 + b^2} \, \| \nabla \xi_{\varepsilon} \|_{L^{\infty}(\Omega)} \right).
\end{equation}
Further explicit computations show that
\begin{equation} \label{norml2q1}
\begin{aligned}
\| \nabla q_{1} \|_{L^{2}(\Omega)} & \leq \| \nabla \xi_{\varepsilon} \|_{L^{2}(\Omega)} + \left\| \dfrac{\partial \xi_{\varepsilon}}{\partial x} \right\|_{L^{2}(\Omega)} + \dfrac{1}{2} \left( \left\| \dfrac{\partial \xi_{\varepsilon}}{\partial y} \right\|_{L^{2}(\Omega)} + \left\| \dfrac{\partial \xi_{\varepsilon}}{\partial z} \right\|_{L^{2}(\Omega)} \right) \\[6pt]
& \hspace{5mm} + \dfrac{(1 + \varepsilon)b}{2} \left( \left\| \nabla \left( \dfrac{\partial \xi_{\varepsilon}}{\partial x} \right) \right\|_{L^{2}(\Omega)} + \left\| \nabla \left( \dfrac{\partial \xi_{\varepsilon}}{\partial y} \right) \right\|_{L^{2}(\Omega)} \right) \\[6pt]
& \hspace{5mm} +  \dfrac{(1 + \varepsilon)c}{2} \left( \left\| \nabla \left( \dfrac{\partial \xi_{\varepsilon}}{\partial x} \right) \right\|_{L^{2}(\Omega)} + \left\| \nabla \left( \dfrac{\partial \xi_{\varepsilon}}{\partial z} \right) \right\|_{L^{2}(\Omega)} \right) \\[6pt]
& \leq \Bigg[ \| \nabla \xi_{\varepsilon} \|_{L^{\infty}(\Omega)} + \left\| \dfrac{\partial \xi_{\varepsilon}}{\partial x} \right\|_{L^{\infty}(\Omega)} + \dfrac{1}{2} \left( \left\| \dfrac{\partial \xi_{\varepsilon}}{\partial y} \right\|_{L^{\infty}(\Omega)} + \left\| \dfrac{\partial \xi_{\varepsilon}}{\partial z} \right\|_{L^{\infty}(\Omega)} \right) \\[6pt]
& \hspace{6mm} + \dfrac{(1 + \varepsilon)b}{2} \left( \left\| \nabla \left( \dfrac{\partial \xi_{\varepsilon}}{\partial x} \right) \right\|_{L^{\infty}(\Omega)} + \left\| \nabla \left( \dfrac{\partial \xi_{\varepsilon}}{\partial y} \right) \right\|_{L^{\infty}(\Omega)} \right) \\[6pt]
& \hspace{6mm} +  \dfrac{(1 + \varepsilon)c}{2} \left( \left\| \nabla \left( \dfrac{\partial \xi_{\varepsilon}}{\partial x} \right) \right\|_{L^{\infty}(\Omega)} + \left\| \nabla \left( \dfrac{\partial \xi_{\varepsilon}}{\partial z} \right) \right\|_{L^{\infty}(\Omega)} \right) \Bigg] \, \sqrt{8 abc \varepsilon^{3}} \, 
\end{aligned}
\end{equation}
and 
\begin{equation} \label{norml2q3}
\begin{aligned}
\| \nabla q_{3} \|_{L^{2}(\Omega)} & \leq \Bigg[ \| \nabla \xi_{\varepsilon} \|_{L^{\infty}(\Omega)} + \left\| \dfrac{\partial \xi_{\varepsilon}}{\partial z} \right\|_{L^{\infty}(\Omega)} + \dfrac{1}{2} \left( \left\| \dfrac{\partial \xi_{\varepsilon}}{\partial x} \right\|_{L^{\infty}(\Omega)} + \left\| \dfrac{\partial \xi_{\varepsilon}}{\partial y} \right\|_{L^{\infty}(\Omega)} \right) \\[6pt]
& \hspace{6mm} + \dfrac{(1 + \varepsilon)a}{2} \left( \left\| \nabla \left( \dfrac{\partial \xi_{\varepsilon}}{\partial x} \right) \right\|_{L^{\infty}(\Omega)} + \left\| \nabla \left( \dfrac{\partial \xi_{\varepsilon}}{\partial z} \right) \right\|_{L^{\infty}(\Omega)} \right) \\[6pt]
& \hspace{6mm} +  \dfrac{(1 + \varepsilon)b}{2} \left( \left\| \nabla \left( \dfrac{\partial \xi_{\varepsilon}}{\partial y} \right) \right\|_{L^{\infty}(\Omega)} + \left\| \nabla \left( \dfrac{\partial \xi_{\varepsilon}}{\partial z} \right) \right\|_{L^{\infty}(\Omega)} \right) \Bigg] \, \sqrt{8 abc \varepsilon^{3}} \, .
\end{aligned}
\end{equation}
The proof is complete after choosing $\varepsilon = \dfrac{L-a}{2a}$ and by inserting the estimates \eqref{estimar2} into \eqref{norml3q1}, \eqref{norml3q3}, \eqref{norml2q1} and \eqref{norml2q3}.
\end{proof}

It is then clear that, by combining Corollary \ref{maincor}, Proposition \ref{newformula} and Proposition \ref{solext}, we are able to derive \textit{explicit} upper bounds on the drag and lift exerted over $K$ by the fluid flow governed by the \text{unique} solution of \eqref{nsstokes0}. We shall give the expression for such bounds in the case where $a=b=c$ and the inflow velocity is given by the following function:
\begin{equation} \label{inflowcos}
h(y,z) = A \cos \left( \dfrac{\pi y}{2L} \right) \cos \left( \dfrac{\pi z}{2L} \right) e_{1} \qquad \forall (y,z) \in \Gamma_{I} \, ,
\end{equation}
for some constant $A > 0$ that represents the magnitude of the inflow.
\begin{figure}[H]
	\begin{center}
		\includegraphics[height=50mm,width=100mm]{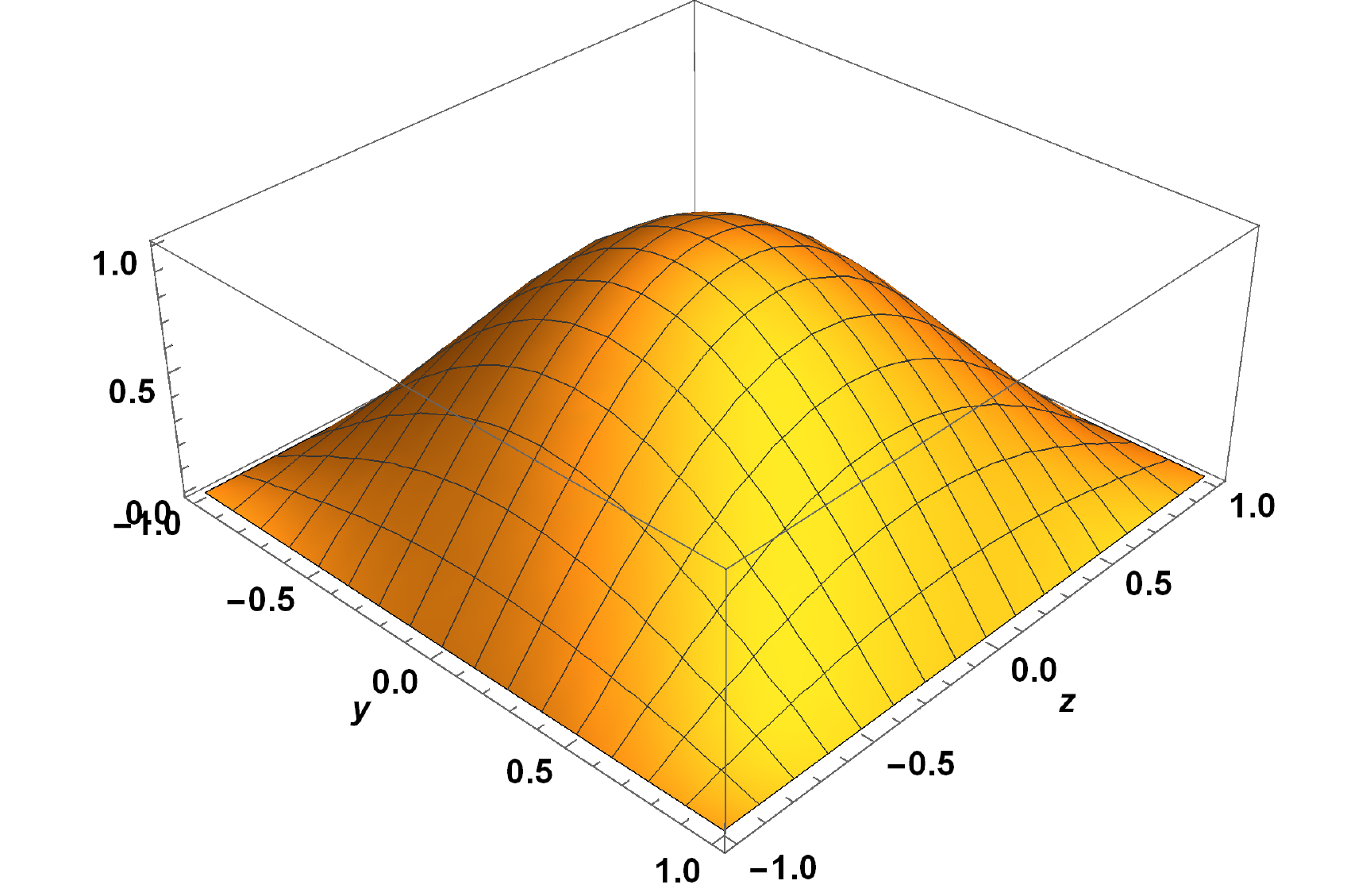} 
	\end{center}
	\vspace*{-3mm}
	\caption{Graph of the function $|h|$ defined in \eqref{inflowcos}, for $A = L = 1$.}\label{inflowgraph}
\end{figure}
\noindent
\begin{corollary} \label{upperdraglift}
Assume \eqref{OmegaR}, define $M>0$ as in \eqref{sigmagamma3d} with $a=b=c$, where 
\begin{equation} \label{conda}
|K| > \dfrac{8}{81} (81 - 16 \pi) L^3 \, .
\end{equation}
Suppose that the boundary datum $h \in H^{1}_{0}(\Gamma_{I})$ is given by \eqref{inflowcos}. If
\begin{equation}\label{finale1}
A \leq \dfrac{\eta^{2}}{4 \sqrt{2} \pi^{1/3}} \, \dfrac{\sqrt{L} \, (L-a)}{\sqrt{2} \, (1+M)L + (L-a)\pi} \, \dfrac{\left[\pi^2 - \dfrac{2}{L^2} \left( \dfrac{3}{2 \pi} (|Q| - |K|) \right)^{2/3} \right]^{1/2}}{\left[ \eta + \dfrac{1}{\sqrt{5} \, \pi} \left( \dfrac{7^{7/3}}{3} \right)^{1/4} \, \left( \dfrac{3}{4} (|Q| - |K|) \right)^{1/6} \right]^{2}} \, ,
\end{equation}
then problem \eqref{nsstokes0} admits a unique weak solution $u \in \mathcal{V}(\Omega)$ such that
\begin{equation} \label{finale11}
\| \nabla u \|_{L^{2}(\Omega)} \leq \sqrt{2L} \, A \left( \dfrac{1+M}{L-a} + \dfrac{\pi}{\sqrt{2} \, L} \right) \left[ \sqrt{2} + \dfrac{1}{\sqrt{5} \, \pi \, \eta} \, \left( \dfrac{7^{7/3}}{3} \right)^{1/4} \, \left( \dfrac{3}{4} (|Q| - |K|) \right)^{1/6} \, \right] .
\end{equation}
Furthermore, the following estimates for the drag and lift applied over $K$ hold:
\begin{equation} \label{finale2}
| \mathcal{D}_{K}(u,p) | \leq \Psi(u) \qquad \text{and} \qquad | \mathcal{L}_{K}(u,p) | \leq \Psi(u) \, ,
\end{equation}
with
\begin{equation} \label{finale3}
\begin{aligned}
\Psi(u) & \doteq 6 \eta \, \sqrt{L-a} \, \left( 2 + \sqrt{3} + \sqrt{82} \, \dfrac{L+a}{L-a} \right) \| \nabla u \|_{L^{2}(\Omega)} \\[6pt] 
& \hspace{4mm} + \left( 1 + \dfrac{3}{4} \sqrt{6} \right) \, \dfrac{ 2 \pi^{1/3} (L+a) }{\left[ \pi^2 - \dfrac{2}{L^2} \left( \dfrac{3}{2 \pi} (|Q| - |K|) \right)^{2/3} \, \right]^{1/2}} \, \| \nabla u \|^{2}_{L^{2}(\Omega)} \, .
\end{aligned}
\end{equation}
\end{corollary}
\noindent
\begin{proof}
First of all, assumption \eqref{conda} ensures that
$$
\min \left\{ \sqrt[3]{\frac{3}{2\pi} (|Q|-|K|)}  \, , \,  \dfrac{4L}{3} \right\} = \sqrt[3]{\frac{3}{2\pi} (|Q|-|K|)} \, ,
$$
see Section \ref{funcine}. Now, the function $h$ in \eqref{inflowcos} is divergence-free, so that in Theorem \ref{teononconstant} we get
$$
\Phi(h) = \sqrt{2L} \, \left( \dfrac{1+M}{L-a} + \dfrac{\pi}{\sqrt{2} \, L} \right) A \, .
$$
In virtue of Corollary \ref{maincor}, condition \eqref{finale1} guarantees the existence of a unique weak solution $u \in \mathcal{V}(\Omega)$ to problem \eqref{nsstokes0} satisfying the estimate \eqref{finale11}. The upper bounds \eqref{finale2}-\eqref{finale3} on the drag and lift follow from Proposition \ref{newformula}, adjusting the expressions given in Proposition \ref{solext} to the case $a=b=c$.
\end{proof}

\bigskip

\par\medskip\noindent
{\bf Acknowledgements.} The Author is supported by the Research Programme PRIMUS/19/SCI/01, by the program GJ19-11707Y of the Czech National Grant Agency GA\v{C}R, and by the University Centre UNCE/SCI/023 of the Charles University in Prague.

\phantomsection
\addcontentsline{toc}{section}{References}
\bibliographystyle{abbrv}
\bibliography{references}
\vspace{5mm}

\begin{minipage}{100mm}
Gianmarco Sperone\\
Department of Mathematical Analysis\\
Faculty of Mathematics and Physics\\
Charles University in Prague\\
Sokolovská 83\\
186 75 Prague - Czech Republic\\
{\bf E-mail}: sperone@karlin.mff.cuni.cz
\end{minipage}
\end{document}